\documentclass{article}
\usepackage{amscd}      % to be able to use "CD" environment
\usepackage{amssymb}     
\usepackage{amsmath}     
\usepackage{xypic}      % to be able to use "diagram" environment of 
\LaTeXdiagrams          % xypic to typeset commutative diagrams 
\usepackage[all,v2]{xy}      
\xyoption{2cell}
\UseAllTwocells
\xyoption{frame}
\CompileMatrices
\allowdisplaybreaks[4] 

\usepackage{theorem}

\newtheorem{prop}{Proposition}[section]  
\newtheorem{lem}[prop]{Lemma}

\newtheorem{cor}[prop]{Corollary}
\newtheorem{them}[prop]{Theorem}

\newtheorem{scholum}[prop]{Scholum}
\newtheorem{nothem}{Theorem}

\theorembodyfont{\upshape}

\newtheorem{defn}[prop]{Definition}

\newtheorem{numrmk}[prop]{Remark}
\newtheorem{numnote}[prop]{Note}
\newtheorem{numex}[prop]{Example}

\newtheorem{note}{Note}

\newtheorem{axiom}[prop]{Axiom}%\renewcommand{\theaxiom}{}

\newtheorem{rmk}{Remark}

\newenvironment{pf}{\begin{trivlist}\item[]{\sc Proof.}}%
            {\nolinebreak $\Box$ \end{trivlist}}

\newcommand{\noprint}[1]{}

\newcommand{\Section}[1]{\newpage\section{#1}}

\renewcommand{\tilde}{\widetilde}

\newcommand{\toto}{\rightrightarrows}
\newcommand{\qed}{{\nolinebreak $\,\,\Box$}}

\newcommand{\upst}{^{\ast}}

\newcommand{\perf}{{\rm pf}}

\newcommand{\lst}{_{\ast}}

\newcommand{\com}{^{\scriptscriptstyle\bullet}}
\newcommand{\lcom}{_{\scriptscriptstyle\bullet}}

\newcommand{\argument}{{{\,\cdot\,}}}

\renewcommand{\AA}{{\mathfrak A}}

\renewcommand{\SS}{{\mathfrak S}}

\newcommand{\RR}{{\mathfrak R}}

\newcommand{\Mm}{{\mathfrak m}}

\newcommand{\Dd}{{\mathfrak d}}

\newcommand{\zz}{{\mathbb Z}}

\newcommand{\aaa}{{\mathbb A}}

\renewcommand{\O}{{\cal O}}

\newcommand{\del}{\partial}
\newcommand{\du}{\Dd}
\newcommand{\resto}{{\,|\,}}
\newcommand{\st}{\mathrel{\mid}}

\newcommand{\gr}{\mathop{\rm gr}\nolimits}

\renewcommand{\div}{\mathop{\rm div}}

\newcommand{\im}{\mathop{\rm im}}
\newcommand{\ob}{\mathop{\rm ob}}

\newcommand{\cok}{\mathop{\rm cok}}

\newcommand{\spec}{\mathop{\rm Spec}\nolimits}

\newcommand{\id}{\mathop{\rm id}\nolimits}
\newcommand{\Hom}{\mathop{\rm Hom}\nolimits}

\newcommand{\shom}{\mathop{\rm Hom}\nolimits^{\scriptscriptstyle\Delta}} 

\newcommand{\Tor}{\mathop{\rm Tor}\nolimits}

\newcommand{\Homu}{\mathop{\underline{\rm Hom}}\nolimits}
\newcommand{\Deru}{\mathop{\underline{\rm Der}}\nolimits}

\newcommand{\injectlim}{\mathop{\lim\limits_{\textstyle\longrightarrow}}\limits}
\newcommand{\projectlim}{\mathop{{\lim\limits_{\textstyle\longleftarrow}}}\limits}
\newcommand{\comp}{\mathbin{{\scriptstyle\circ}}}
\newcommand{\ol}{\overline}

\newcommand{\longiso}{\stackrel{\textstyle\sim}{\longrightarrow}}

\newcommand{\doublearrowstack}[2]%
                      {{{{\scriptstyle#1}\atop{\textstyle\longrightarrow}}\atop{{\textstyle\longrightarrow}\atop{\scriptstyle#2}}}}
\newcommand{\rightleftarrowstack}[2]%
                      {{{{\scriptstyle#1}\atop{\textstyle\longrightarrow}}\atop{{\textstyle\longleftarrow}\atop{\scriptstyle#2}}}}
\newcommand{\leftrightarrowstack}[2]%
                      {{{{\scriptstyle#1}\atop{\textstyle\longleftarrow}}\atop{{\textstyle\longrightarrow}\atop{\scriptstyle#2}}}}

\newcommand{\ses}[5]%
{0\longrightarrow#1\stackrel{#2}{ \longrightarrow}#3\stackrel{#4}{
\longrightarrow}#5\longrightarrow0} 

\newcommand{\dt}[6]%
{#1\stackrel{#2}{longrightarrow}#3 \stackrel{#4}{\longrightarrow}#5
\stackrel{#6}{\longrightarrow} #1[1]}  
 
\newcommand{\cat}[1]%
{(\mbox{\rm #1})}

\title{{\bf \LARGE Differential Graded Schemes~I:}  \\
{\bf \Large  Perfect Resolving Algebras}}
\author{K. Behrend} 
\date{December 16, 2002}

\begin{document}
\sloppy
\maketitle

\begin{abstract}
We introduce {\em perfect resolving algebras }and study their fundamental
properties.  These algebras are basic for our theory of 
differential graded schemes, as they give rise to {\em affine }differential
graded schemes.  We also introduce {\em \'etale morphisms}.  The
purpose for studying these, is that they will be  used
to glue differential graded schemes from affine ones with respect to
an \'etale topology.
\end{abstract}

\tableofcontents

\newpage
\section*{Introduction}
\addcontentsline{toc}{section}{Introduction}

This is the first in a series of papers devoted to establishing a
workable theory of differential graded schemes.

Here we lay the necessary algebraic foundations for this theory.
Differential graded schemes will be glued with respect to an \'etale
topology from affine differential graded schemes, very much like usual
schemes are glued with respect to the Zariski topology from usual
affine schemes, or algebraic spaces are glued with respect to the
\'etale topology from affine schemes.

Thus our goal in this paper is twofold: we introduce an
appropriate class of differential graded algebras providing us with a
good class of affine differential graded schemes, and we 
introduce the notion of \'etale morphism between such differential
graded algebras.

The first principle we follow, is that all differential graded
algebras which represent geometric objects, are 

(i) graded commutative with 1, over a field of characteristic 0,

(ii) graded {\em in
non-positive degrees}, if the differential has degree $+1$, which is
the convention we will follow.  

A consequence of this is that
every such differential graded algebra $A$ has a morphism of
differential graded algebras $A\to h^0(A)$, where $h^0(A)$ is the 0-th
cohomology module of $A$.

A {\em quasi-isomorphism }of differential graded algebras is a morphism
$A\to B$ which induces an isomorphism on cohomology modules
$h^i(A)\to h^i(B)$, for all $i\leq0$. 

Another principle is, that quasi-isomorphic differential graded
algebras should give rise to identical geometric objects.  Thus we may
replace the arbitrary differential graded algebra $A$ (concentrated in
non-positive degree) by a quasi-isomorphic differential graded algebra
$A'$, which is free as a graded algebra, disregarding its differential
(a property which has been referred to as {\em quasi-free }in the
literature). (See Scholum~\ref{skol} for a proof of the
existence of {\em resolutions }$A'\to A$.)

Thus, we are able to restrict our attention do differential graded
algebras which are free as graded commutative algebras with 1, free on
a set of generators all of which have non-positive degree.  We call
algebras satisfying this property {\em resolving algebras}, because
their purpose is to resolve more general differential graded algebras
(see Definition~\ref{def.res}).
We also required a term which is shorter than `quasi-free on a set of
generators in non-positive degree', and can more easily be qualified.

For purposes of  our geometric theory, we need  finiteness assumptions
on resolving algebras.  Thus we call a resolving algebra {\em finite
}if we can find a finite set of quasi-free generators for it. The
finite resolving algebras will serve as a category of local models for
differential graded schemes.  In other words, every differential
graded scheme will be locally (with respect to the \'etale topology)
given by a finite resolving algebra.  Put another way, every
differential graded scheme will locally determine up to
quasi-isomorphism a finite resolving algebra as an analogue of `affine
coordinate ring'.  

On the other hand, it turns out that (the differential graded schemes
associated to) finite resolving algebras form too small a class to
be considered as {\em  affine }differential graded scheme. 
More specifically, a fundamental (and too useful to forgo) property of
affine morphisms of usual schemes is, that the affine property is {\em
local }in the base, or target, of the morphism.  (If local means local
with respect to the \'etale topology, this is one of the most
important results of \'etale descent theory.)

Thus we are led to relax the requirement of finiteness on resolving
algebras.  We call a resolving algebra $A$ {\em perfect }(see
Definition~\ref{def.per}) if 

(i) it is {\em quasi-finite}, i.e., we can find a set of quasi-free
generators, with finitely many elements of every degree,

(ii) its differential graded module of differentials $\Omega_A$ (see
Section~\ref{difcot}) tensored with $h^0(A)$ is a perfect complex of
$h^0(A)$-modules, i.e., is Zariski locally in $\spec h^0(A)$
quasi-isomorphic to a finite complex of finite rank free modules.

Note that every finite resolving algebras is perfect. 

Using  perfect resolving algebras to define affine differential
graded schemes, the notion of affine morphism of differential
graded schemes is local in the \'etale topology on the base (see
\cite{dgsII}).  Thus we choose perfect resolving algebras as our
affine models for differential graded schemes.

We prove two fundamental facts about perfect resolving algebras,
showing that they are, in fact, not very far from finite resolving
algebras:

(i) Every perfect resolving algebra $A$ is locally with respect to the
Zariski topology on $\spec h^0(A)$ quasi-isomorphic to a finite
resolving algebra (see Theorem~\ref{loc.fin}).

(ii) The derivations of a  perfect resolving algebra are in a certain
sense compatible with the derivations of its truncations (the
subalgebras generated by finite subsets of a generating set). See
Theorem~\ref{plim} for the precise statement.  A similar result also
holds for homotopy groups (see Corollary~\ref{greatc}). 

\subsubsection{The main results}

There is a natural structure of simplicial closed model category on
the category of differential graded algebras.  Resolving algebras are
cofibrant objects for this closed model category structure and 
for any two resolving algebras $A$, $B$, the simplicial set of
morphisms from $B$ to $A$, denoted $\shom(B,A)$, is fibrant, i.e., has
the Kan property, and can thus be considered as a (topological)
space. In particular, we have homotopy groups $\pi_\ell\shom(B,A)$,
for $\ell>0$.  These facts are reviewed in~\ref{sec.scmc}.
In~\ref{sec.homo}, we prove that a morphism of resolving algebras is a
quasi-isomorphism if and only if it is a homotopy equivalence. 

Given two resolving algebras $A$, $B$, we also have the differential graded
$A$-module $\Deru(B,A)$: an element $D:B\to A$ of degree $n$ in
$\Deru(B,A)$ is a degree $n$ homomorphism of graded vector spaces,
which satisfies the graded Leibniz rule (see
Definition~\ref{def.imd}).  

The main results of this paper relate the homotopy groups
$\pi_\ell\shom(B,A)$ to the $h^0(A)$-modules $h^{-\ell}\Deru(B,A)$:

\begin{nothem}
Let $A$ and $B$ be perfect resolving algebras. Then, for every
$\ell>0$,  there is a canonical bijection
$$\Xi_\ell:h^{-\ell}\Deru(B,A)\longrightarrow \pi_\ell\shom(B,A)\,.$$
The bijection $\Xi_\ell$ is an isomorphism of groups in the following
two cases:

(i) if $\ell\geq2$,

(ii) if $B$ is generated by a set of homogeneous generators, all of
which have the same degree.
\end{nothem}

There is also a relative version of this theorem, making Case~(ii)
more interesting and useful. In fact, the relative version and
Case~(ii) are the key results, as they allow proofs by induction. The
relative version of Case~(ii) also 
extends to $\ell=0$ (see Corollary~\ref{soepf}). 

The map $\Xi_\ell$ is defined by Formula~(\ref{def.xil})
in~\ref{sec.lin}.  The theorem is Theorem~\ref{Bij} and its
Corollary~\ref{greatc}. 

\subsubsection{Overview}

In Section~\ref{sec.DGA} we start by reviewing basic definitions
involving differential graded algebras and differential graded
modules.  We introduce resolving algebras. We also study derivations,
differentials and the cotangent complex.

In Section~\ref{sec.etale} we introduce the notion of \'etale morphism
between quasi-finite resolving algebras (see
Definition~\ref{defetale}).  There are many different ways to
characterize \'etale morphisms.  The most important are as
follows. A morphism $A\to B$ is \'etale, if and only if any of the
following equivalent conditions holds:

(i) the relative cotangent complex $L_{B/A}$ is acyclic,

(ii) $\spec h^0(B)\to\spec h^0(A)$ is an \'etale morphism of usual
affine schemes 
and $h^i(A)\otimes_{h^0(A)}h^0(B)\to h^i(B)$ is an isomorphism, for
all $i\leq0$, (in particular, any quasi-isomorphism is \'etale),

(iii) (if $A$ and $B$ are perfect) the induced morphism of
completions $\widehat A\to \widehat B$
is a quasi-isomorphism, for the completions at every
augmentation of $B$
(see Theorem~\ref{et.com.crit}), 

(iv) (if $A$ and $B$ are perfect) for every $B\to C$, the map
$\shom(B,C)\to\shom(A,C)$ induces isomorphisms on homotopy groups
$\pi_\ell$, for all, or for one fixed $\ell>0$ (see
Proposition~\ref{etahom}). 

The most important example of an \'etale morphism is the {\em standard
}\'etale morphism $A\to B$, where $B$ is generated over $A$ by formal 
variables $x_1,\ldots,x_r$ in degree 0 and $\xi_1,\ldots,\xi_r$ in
degree -1. The differential on $B$ is given by $dx_i=0$ and
$d\xi_j=f_j\in A^0[x_1,\ldots,x_r]$. Moreover, 
$\det(\frac{\del f_j}{\del x_i})$ is a unit in
$h^0(A)[x_1,\ldots,x_r]/(f_1,\ldots,f_r)$. 

We prove that every \'etale morphism is locally standard (see
Proposition~\ref{loc.str}, for the exact statement). Thus the study of
\'etale morphisms can often be reduced to the study of standard
\'etale morphisms.

As a byproduct of the equivalence of~(i) and~(ii) we get a very useful
quasi-isomorphism criterion: $A\to B$ is a quasi-isomorphism if and
only if $L_{B/A}$ is acyclic and $h^0(A)\to h^0(B)$ an
isomorphism. (See Corollary~\ref{qiscondition}.)

In Section~\ref{sec.perfect}, we introduce perfect resolving algebras
and prove the two fundamental facts alluded to, above.  

Section~\ref{sec.linear}, is devoted to the  proof of the main theorem
mentioned above, 
to the effect that we can `linearize' homotopy groups. 

Finally, in Section~\ref{sec.fibered}, we prove that any morphism
$A\to B$ between finite resolving algebras admits a factorization
$A\to B\to B'$, where $A\to B$ makes $B$ a finite resolving algebra
over $A$ and $B'\to B$ is a quasi-isomorphism. As an application, we
prove that the fibered homotopy coproduct (or derived tensor product,
as we prefer to call it) of a diagram of finite (perfect) resolving
algebras can be represented, again, by a finite (perfect) resolving
algebra.

\subsubsection{Acknowledgements}

I wish to thank the Research Institute for the Mathematical Sciences
in Kyoto, for providing a wonderful work environment during my visit
in 1999/2000. I would especially like to thank my host, Professor
K. Saito, for the warm hospitality.  Most of the work on \'etale
morphisms was done while I was at RIMS for a sabbatical in 1999/2000.

I would also like to thank the Mathematical Sciences Research
Institute in Berkeley. During a stay there in early 2002, I discovered
the importance of perfect resolving algebras and worked out the
results on linearization of homotopy groups.

This work was also continuously supported by a Research Grant from the
National Science and Engineering Council of Canada. 

Finally, I would like to thank I.~Ciocan-Fontanine, E.~Getzler,
M.~Kapranov, B.~Toen and G.~Vezzosi for helpful discussions.

\Section{Differential Graded Algebras}\label{sec.DGA}

\subsection{Review and terminology}
 
We will fix a base field of characteristic zero, denoted $k$,
throughout.  All rings and algebras are assumed to be commutative with
unit (if they are graded, then commutative means graded
commutative). Algebra without qualifier means $k$-algebra. If $a$ is a
homogeneous element of a graded $k$-vector space, we denote by
$\bar{a}$ its degree.

\paragraph{Differential graded algebras}

A {\em differential graded }$k$-algebra $A$ is a graded $k$-algebra 
$$A^\natural=\bigoplus_{n\in\zz}A^n$$ endowed with a differential
$d:A^\natural\to A^\natural$ of degree 1, i.e., a degree 1
homomorphism of graded $k$-vector spaces satisfying $d^2=0$ and the
graded Leibniz rule $$d(ab)=(da)b+(-1)^{\bar{a}}a\, db.$$ We always
use the notation $A^\natural$ for the underlying graded algebra of the
differential graded algebra $A=(A^\natural,d)$.

For a differential graded algebra $A$, the cohomology $h\upst(A)=\ker
d/\im d$ is a graded $k$-algebra. The degree 0 cohomology $h^0(A)$ is
a subalgebra.  If $A$ is concentrated in non-positive degrees, then
there is a morphism of differential graded algebras $A\to h^0(A)$.

\paragraph{Modules}

A {\em differential graded module }$M$ over the differential graded
algebra $A$ is a (left) graded $A^\natural$-module $M^\natural$
endowed with a differential $d_M:M^\natural\to M^\natural$ of degree
1, i.e., a degree 1 homomorphism of graded $k$-vector spaces
satisfying $d_M^2=0$ and the graded Leibniz rule
$$d_M(am)=(da)m+(-1)^{\bar{a}}a\, d_M m,$$ for $a\in A$ and $m\in
M$. Note that $h\upst(M)$ is a graded $h\upst(A)$-module.

Sometimes, we will write the action of $A$ on the right. In this case
the notation $ma$ is understood to mean
\begin{equation}\label{understood}
ma:=(-1)^{\bar{a}\bar{m}}am\,.
\end{equation}

Note that every differential graded $A$-module has an underlying
complex of $k$-vector spaces.

Let $M$ and $N$ be differential graded $A$-modules. Tensor product and
internal hom are defined as follows.  The {\em tensor product }$M\otimes_A
N$ has underlying graded $A^\natural$-module
$$M^\natural\otimes_{A^\natural}N^\natural$$
and the differential is given by 
$$d_{M\otimes N}(m\otimes n)=d_Mm\otimes n + (-1)^{\bar{m}}m\otimes
d_Nn.$$
The {\em internal hom }$\Homu_A(M,N)$ has underlying graded
$A^\natural$-module
$$\Homu_A(M,N)^\natural=\bigoplus_n
\Hom^n_{A^\natural}(M^\natural,N^\natural)$$ 
and differential given by 
$$d_N\phi(m)=d_{\Homu(M,N)}(\phi)(m)+(-1)^{\bar{\phi}}\phi(d_Mm),$$
for $m\in M$ and $\phi\in\Homu_A(M,N)$. The $k$-vector space of
differential graded $A$-module homomorphisms $M\to N$ is the set of
0-cocycles in $\Homu_A(M,N)$:
$$\Hom_A(M,N)=Z^0\Homu_A(M,N)$$

Note that differential graded $A$-modules form an abelian
category with kernels, cokernels, images and direct
sums taken degree-wise.  In fact, the category of differential graded
$A$-modules is an abelian subcategory of the category of complexes of
$k$-vector spaces.

We have the following formulas:
$$\Homu_A(M,N)\otimes_AP=\Homu_A(M,N\otimes_AP)\,,$$
$$\Homu_A(M\otimes_AN,P)=\Homu_A\big(M,\Homu_A(N,P)\big)\,,$$
for differential graded modules $M$, $N$, $P$ over the differential
graded algebra $A$,
$$\Homu_A(M,N)=\Homu_B(M\otimes_AB,N)\,,$$
for a morphism of differential graded algebras $A\to B$ and a
differential graded modules $M$ over $A$ and $N$ over $B$.

\paragraph{Cones}

Given a differential graded $A$-module $M$, the shift $M[1]$ is
defined by shifting the underlying complex of $k$-vector spaces:
$$(M[1])^i=M^{i+1};\quad d_{M[1]}=-d_M.$$
The shift is again a differential graded $A$-module in a natural way. 

The cone over the homomorphism $\phi:M\to N$ of differential graded
$A$-modules is defined by 
$$C(\phi)=C(M\to N)=N\oplus M[1]. $$
Thus $C(M\to N)$ is again a differential graded $A$-module. There is a
canonical triangle
$$M\stackrel{\phi}{\longrightarrow}N\longrightarrow
C(\phi)\longrightarrow M[1]$$
of differential graded $A$-modules, which induces a long exact sequence
of cohomology groups.

Given two homomorphisms of differential graded $A$-modules 
$$M\stackrel{\phi}{\longrightarrow}N\stackrel{\rho}{\longrightarrow}
P,$$
then, if $\rho\comp\phi=0$, there is a canonical homomorphism
$C(\phi)\to P$ (induced by the second projection) making
$$\begin{diagram}
N\rto\drto_\rho & C(\phi)\dto\\
& P\end{diagram}$$
commute.  If 
$$0\longrightarrow M\stackrel{\phi}{\longrightarrow} N\longrightarrow
P\longrightarrow 0$$ is exact, then $C(\phi)\to P$ is a
quasi-isomorphism.

\paragraph{Spectral Sequence}

Let $M$ be a differential graded $A$-module.  Then $M$ is filtered by
differential graded submodules $F^pM= A\cdot M^{\geq p}$. Thus we get an
associated spectral sequence, which converges if both $A$ and $M$ are
bounded above. 

\begin{prop}\label{spec.se}
Assume the differential graded algebra $A$ is concentrated in
non-positive degrees.  Let $M$ be a differential graded $A$-module
which is bounded above, and such that $M^\natural$ is a free
$A^\natural$-module. Then we have a spectral sequence of
$h^0(A)$-modules
$$E_2^{p,q}=h^q(A)\otimes_{h^0(A)}h^p\big(M\otimes_A h^0(A)\big)
\Longrightarrow h^{p+q}(M)\,.\quad\Box$$
\end{prop}

As an application, we may for example deduce, that $M$ is acyclic if
and only if $M\otimes_A h^0(A)$ is acyclic.

\subsection{Resolving algebras}

\paragraph{Symmetric algebras}

Given a complex of $k$-vector spaces $V$, the $n$-th symmetric power,
notation $S^nV$, is defined to be the quotient of $V^{\otimes n}$ by
all relations of the form 
$$x_1\otimes\ldots\otimes x_p\otimes x\otimes y\otimes
y_1\otimes\ldots\otimes y_q=(-1)^{\ol{x}\ol{y}}x_1\otimes\ldots\otimes
x_p \otimes 
y\otimes x\otimes y_1\otimes\ldots\otimes y_q\,,$$
for $x$, $y$, $x_1,\ldots,x_p$, $y_1\ldots,y_q$ homogeneous elements
of $V$ and $p+q+2=n$. Thus $S^nV$ is again a
complex of $k$-vector spaces, for all $n\geq0$.  The direct sum
$SV=\bigoplus_{n\geq0}S^nV$ is a differential graded algebra. The
functor $V\mapsto SV$ is a left adjoint for the forgetful functor
$$(\text{differential graded algebras})\longrightarrow(\text{complexes
of $k$-vector spaces})\,.$$

\begin{defn}\label{free.alg}
A differential graded algebra $A$ which is isomorphic to $SV$, for
some complex of $k$-vector spaces $V$, is called {\bf free}. If
$V\subset A$ is a subcomplex inducing an isomorphism $SV\to A$, we
call $V$ a {\bf free basic complex }for $A$.

If $A\to B$ is a morphism of differential graded algebras and there
exists a subcomplex $V\subset B$ such that $A\otimes_k SV\to B$ is an
isomorphism, then we call $A\to B$ {\bf free}, or we say that $B$ is
free over $A$.  Moreover, $V$ is called a {\bf free basic complex }for
$B$ over $A$.
\end{defn}

The importance of free differential graded algebras for us is that
they occur as tangent spaces of differential graded schemes. 

\paragraph{Quasi-free algebras}

As a special case of symmetric algebras, we may consider symmetric
algebras $SV$, on a graded $k$-vector spaces $V$. This is the case of
complexes $V$ with zero differential.  In this case $SV$ is simply a
graded $k$-algebra, as it has vanishing differential, too.  The
functor $V\mapsto SV$ is a left adjoint for the forgetful functor
$$(\text{graded $k$-algebras})\to(\text{graded $k$-vector
spaces})\,.$$ 
Another common notation for $SV$ is $k[V]$.  If $(x_i)$ is a
homogeneous basis for the graded $k$-vector space $V$, then we denote
$k[V]$ also by $k[x]$.

\begin{defn}
Let $A$ be a differential graded algebra. Suppose that $V\subset
A$ is a graded sub-$k$-vector space, such that $SV\to A^\natural$ is
an isomorphism of graded $k$-algebras. Then we say that the
differential graded algebra $A$ is {\bf quasi-free }and we call $V$ a {\bf
basic }space for $A$.

Let $V\subset A$ be a basic space for the quasi-free
differential graded algebra $A$.  If $(x_i)$ is a homogeneous
$k$-basis for the graded $k$-vector space $V$, then we call $(x_i)$ a
{\bf basis }for $A$.

Let $A\to B$ be a morphism of differential graded algebras. If
$V\subset B$ is a graded subspace such that $A^\natural\otimes_k SV\to
B^\natural$ is an isomorphism, then $B$ is {\bf quasi-free }over $A$
and $V$ is a {\bf basic }space for $B$ over $A$.  If $(x_i)$
is a homogeneous basis for $V$, then it is called a {\bf basis
}for $B$ over $A$.
\end{defn}

Thus any basis $(x_i)$ of a quasi-free differential graded algebra $A$
defines an isomorphism $k[x]\to A^\natural$.

We have chosen the terms {\em basic }and {\em basis}, rather than the
more logical terms {\em quasi-basic }and {\em quasi-basis}, because
these terms will be used much more often than the terms defined in
Definition~\ref{free.alg}, and we do not want the prefix `quasi' to
take over the paper.

Let $B$ be quasi-free over $A$ and $(x_i)$ a basis. If $A\to C$
is a morphism of differential graded algebras, then any family $(c_i)$
of homogeneous elements of $C$, such that $\deg c_i=\deg x_i$ for all
$i$, induces a unique morphism of graded algebras $B^\natural
\to C^\natural$, extending $A^\natural\to C^\natural$.  This morphism
$B\to C$ is a morphism of differential graded algebras if and only if
it maps $dx_i$ to $dc_i$ for all $i$.

The notion of quasi-freeness itself is not very useful for us.
Its main purpose is to enable the definition of {\em resolving
algebra}, which is next.

\paragraph{Resolving algebras}

We now come to be most important set of concepts for this work.

\begin{defn}\label{def.res}
We call a differential graded algebra $A$ a {\bf resolving algebra},
if it is quasi-free and there exists a basis $(x_i)$ for $A$,
such that $\deg x_i\leq0$, for all $i$.

A morphism of differential graded algebras $A\to B$ is called a {\bf
resolving morphism}, if $B$ is quasi-free over $A$ and  there exists a
basis $(x_i)$ for $B$ over $A$, such that $\deg x_i\leq0$, for
all $i$. 

If we speak of a {\bf basis }for a resolving algebra or a resolving
morphism, it is understood that it consists of elements $x_i$ such
that $\deg x_i\leq0$, for all $i$.
\end{defn}

The main purpose of resolving morphisms for us is that they are {\em
cofibrations }for the natural simplicial closed model category
structure on the category of differential graded algebras.  For geometric
purposes, we have to put some additional finiteness assumptions:

\begin{defn}
A resolving algebra $A$ is called {\bf quasi-finite}, if there exists a
basis $(x_i)$ for $A$ satisfying

(i) for every $n>0$, the set $\{i\st\deg x_i=n\}$ is empty,

(ii) for every $n\leq0$, the set $\{i\st\deg x_i=n\}$ is finite.

\noindent Any basis for $A$ satisfying these two properties is
called a {\bf quasi-finite }basis or a {\bf coordinate system }for
$A$.

A resolving morphism $A\to B$ is called {\bf quasi-finite}, if there
exists a {\bf quasi-finite }basis for $B$ over $A$, i.e., a
basis  $(x_i)$, satisfying~(i) and~(ii).
\end{defn}

\begin{defn}
A resolving algebra $A$ is called {\bf finite}, if there exists a
finite basis $(x_i)$ for $A$ such that $\deg x_i\leq0$, for all $i$.
Any such basis for $A$ is called a {\bf finite }basis, or a {\bf
finite coordinate system}. 

A resolving morphism $A\to B$ is called {\bf finite}, if there exists
a {\bf finite }basis for $B$ over $A$, i.e., a basis $(x_i)$ which is
finite and satisfies $\deg x_i\leq 0$, for all $i$.
\end{defn}

The purpose of finite resolving algebras for us is, that they provide
local models for differential graded schemes.  There is another
important class of resolving algebras which are somewhere between
finite and quasi-finite resolving algebras.  These are the {\em
perfect }resolving algebras introduced in
Definition~\ref{def.per}. Perfect resolving
algebras serve as affine differential graded schemes.

\begin{defn}
Let $A$ be a finite resolving algebra.  If we can find a finite basis
$(x_i)$ for $A$ such that $\deg x_i\geq -n$, for all i, then we
say that $A$ is of {\bf amplitude }$n$.

If $A\to B$ is a finite resolving morphism and we can find a finite
basis $(x_i)$ for $B$ over $A$ such that $\deg x_i\in[-n,0]$, for all
$i$, then we say that $A\to B$ has {\bf amplitude }$n$. 
\end{defn}

The reason for the terminology is that resolving algebras serve as
resolutions:

\begin{defn}
If $C$ is a differential graded algebra and $A\to C$ is a morphism of
differential graded algebras, where $A$ is a resolving algebra, then
we call $A\to C$ a {\bf resolution }of $C$, if $A\to C$ is a
quasi-isomorphism. 

If $A\to C$ is a morphism of differential graded algebras, $A\to B$ a
resolving morphism, and $B\to C$ a quasi-isomorphism over $A$, then
$A\to B\to C$ is a {\bf resolution }of $A\to C$.
\end{defn}

Let us show that morphisms between resolving algebras always admit 
resolutions:

\begin{prop}\label{res.ex}
For any morphism of resolving  algebras $A\to C$,  there exists
a resolution $A\to B\to C$.  If $A$ and $C$ are quasi-finite, we may
choose $A\to B$ quasi-finite.
\end{prop}
\begin{pf}
We  construct inductively a sequence of partial resolutions $A\to
B_{(n)}\to C$ with the properties:

(i) $A\to B_{(n)}$ is resolving and the composition $A\to B_{(n)}\to
C$ is equal to $A\to C$,

(ii) $h^{-n}(B_{(n)})\to h^{-n}(C)$ is surjective,

(ii) $h^{1-n}(B_{(n)})\to h^{1-n}(C)$ is bijective.

\noindent We may take $A$ for $B_{(-1)}$. Once we have constructed
$B_{(n)}\to C$, we choose elements $b_i\in Z^{-n}(B_{(n)})$ generating
the kernel of $h^{-n}(B_{(n)})\to h^{-n}(C)$, and $e_i\in C^{-n-1}$ such that
$de_i$ is the image of $b_i$ under $B_{(n)}\to C$, for all $i$. We
also choose  $c_j\in Z^{-n-1}(C)$ 
generating $h^{-n-1}(C)$. Then
we define $B_{(n+1)}$ by adjoining to 
$B_{(n)}$ formal variables $x_i$ and $y_j$ of degree $-n-1$ and
setting $dx_i=b_i$ and $dy_j=0$.  We define the morphism $B_{(n+1)}\to
C$ by extending $B_{(n)}\to C$ by  $x_i\mapsto
e_i$ and $y_j\mapsto c_j$.  

Finally, we let $B=\injectlim B_{(n)}$. 

To deal with the quasi-finite case, let us make the remark that if $A$
is a quasi-finite resolving algebra, then $h^0(A)$ is a finite type
$k$-algebra and $h^i(A)$ is a finitely generated $h^0(A)$-module, for
all $i$.  Now we examine the above proof more closely, and specify more
carefully what we mean by `generating'. In fact, to
construct $B_{(0)}$ we choose generators for $h^0(C)$ as an
$h^0(A)$-algebra. To construct all $B_{(n+1)}$ for $n\geq0$, we choose
generators of $\ker\big(h^{-n}(B_{(n)})\to h^{-n}(C)\big)$ and
$h^{-n-1}(C)$ as $h^0(B_{(n)})$-modules. Since $h^0(B_{(n)})\to h^0(C)$
is onto if $n\geq0$, we see that in all cases we have finite
generation.
\end{pf}

\begin{scholum}\label{skol}
If $A$ is a differential graded algebra concentrated in non-positive
degrees, such that $h^0(A)$ is a $k$-algebra of finite type and
$h^i(A)$ is a 
finitely generated $h^0(A)$-module, for every $i$, then there exists a
resolution $A'\to A$, where $A'$ is quasi-finite.
\end{scholum}

\subsection{The simplicial closed model category
structure}\label{sec.scmc}

For the definitions and properties of simplicial sets and closed model
categories, the reader may consult, for example, the recent text book
by Goerss-Jardine~\cite{simhomthe}.  We will commit the common abuse
of calling a simplicial set a {\em space}. 

\begin{defn}
A {\bf simplicial category }is a category enriched over simplicial
sets.  Thus a simplicial
category $\SS$ is given by

(i) a class of objects $\ob\SS$,

(ii) for any two objects $U$, $V$ of $\SS$ a simplicial set
$\shom(U,V)$,

(iii) for any three objects $U$, $V$, $W$ of $\SS$ a simplicial map
\begin{align*}
\comp:\shom(V,W)\times\shom(U,V)& \longrightarrow \shom(U,W)\\*
(f,g)&\longmapsto f\comp g
\end{align*}

(iv) for every object $U$ of $\SS$, a 0-simplex $\id_U$ in
$\shom(U,U)$, which we can also view as a simplicial map
\begin{equation*}
\ast \stackrel{\id_U}{\longrightarrow} \shom(U,U),
\end{equation*}
such that

(i) the composition $\comp$ is associative, i.e., for four objects
$U$, $V$, $W$, $Z$ of $\SS$ the diagram of simplicial sets
$$\begin{diagram}
{\shom(W,Z)\times\shom(W,V)\times\shom(V,U)}\dto\rto &
{\shom(V,Z)\times\shom(U,V)}\dto\\
{\shom(W,Z)\times\shom(U,W)}\rto & {\shom(U,Z)}
\end{diagram}$$
commutes,

(ii) the objects $\id_U$ act as identities for $\comp$, i.e.,
the diagrams 
$$\xymatrix@C=1pc{
{\shom(U,V)}\dto\drto & {\shom(U,V)}\drto\rto &
{\shom(V,V)\times\shom(U,V)}\dto\\
{\shom(U,V)\times\shom(U,U)}\rto& {\shom(U,V)} & {\shom(U,V)}}$$
commute.
\end{defn}

Passing from $\shom(U,V)$ to the set of 0-simplices $\shom_0(U,V)$ we
get the underlying category $\SS_0$ of the simplicial category $\SS$.

Assume that the underlying category of the simplicial category $\SS$ has a
closed model category structure.  The following property (which $\SS$
might enjoy, or not) is called the {\em simplicial model category
axiom}

\begin{axiom}\label{ax}
 If $j:A\to B$ is a cofibration and $q:X\to Y$ a
fibration then 
$$\shom(B,X)\longrightarrow\shom(A,X)\times_{\shom(A,Y)}\shom(B,Y)$$
is a fibration of simplicial sets, which is trivial if $j$ or $q$ is
trivial.
\end{axiom}

\begin{defn}
If the underlying category of the simplicial category $\SS$ is a
closed model category and the simplicial model category axiom is
satisfied, then we call $\SS$ a {\bf simplicial closed model
category}.  (Note that this notion is weaker than the one treated in
\cite{simhomthe}.)
\end{defn}

Differential graded algebras form a simplicial closed model category:
Let $\AA$ be the category of all differential graded $k$-algebras.

\begin{prop}\label{dgascmc}
Call a morphism $f:A\to B$ in $\AA$ a {\em\bf fibration }if $f$ is
degree-wise surjective.  Call $f$ a {\em\bf weak equivalence }if it is
a quasi-isomorphism, i.e., if it induces bijections on cohomology
groups.  Finally, call $f$ a {\em\bf cofibration }if it satisfies the
{left lifting property }with respect to all trivial fibrations.
With these definitions $\AA$ is a closed model category.

Let $\Omega_n=\Omega(\Delta^n)$ be the algebraic de~Rham complex of
the algebraic 
$n$-simplex
$$\spec k[x_0,\ldots,x_n]/{\textstyle\sum} x_i=1,$$
(which is a differential graded $k$-algebra).
For two differential graded algebras $A$, $B$ define the simplicial
set $\shom(A,B)$ to have the $n$-simplices
$$\shom_n(A,B)=\Hom(A,B\otimes\Omega_n).$$
With this definition, $\AA$ is a simplicial closed model category.
\end{prop}
\begin{pf}
This (and much more) is proved in \cite{hinich}.
\end{pf}

\begin{numrmk}\label{fin.tensor}
By definition, we have
$$\Hom\big(\Delta^n,\shom(A,B)\big)=
\Hom\big(A,B\otimes\Omega(\Delta^n)\big)\,.$$ 
More generally, for every finite simplicial set $K$ we have an
algebraic de Rham complex $\Omega(K)$ and there is a natural bijection
$$\Hom\big(K,\shom(A,B)\big)=\Hom\big(A,B\otimes\Omega(K)\big)\,.$$
For the definition of $\Omega(K)$, see~\cite{BoG}.
\end{numrmk}

We shall now identify a class of cofibrations in $\AA$.

\begin{prop}\label{res.cof}
Any resolving morphism $\iota:A\to B$ of differential graded algebras
is a cofibration.  
\end{prop}
\begin{pf}
We have to show that
$\iota:A\to B$ satisfies the left lifting property with respect to all
surjective quasi-isomorphisms. So let $\pi:C\to D$ be a surjective
quasi-isomorphism of differential graded algebras. Assume given the
commutative diagram of solid arrows
$$\begin{diagram}
{A\dto_\iota \rto^g} & {C\dto^\pi}\\
{B\rto_f\urdotted|>\tip^h} & D.
\end{diagram}$$
We need to show the existence of the dotted arrow $h$. 

Let $(y_i)_{i\in I}$ be a basis for $B$ over $A$. By induction, we may
assume that all $y_i$, $i\in I$, 
have the same degree, say $n\leq0$. This implies that $dy_i\in A$, for
all $i\in I$, or more precisely, that there exist $a_i\in A$ such that
$\iota a_i=d y_i$. Moreover, $da_i=0$. (Note that this is where we use
the assumption of non-positive degree.)

Since $C\com\to D\com$ is a quasi-isomorphism, the diagram
$$\begin{diagram}
{C^n/dC^{n-1}\dto\rto^d} & {Z^{n+1}(C\com)\dto}\\
{D^n/dD^{n-1}\rto^d} & {Z^{n+1}(D\com)}
\end{diagram}$$
of $k$-vector spaces is cartesian. So by surjectivity of
$\pi:C^{n-1}\to D^{n-1}$ we have that $C^n$ maps onto the fibered
product
$$\begin{diagram}
{C^n\rto|>>\tip} &
{{\phantom{d}\cdot\phantom{d}}\dto_\pi\rto^{d\phantom{MM}}} & 
{Z^{n+1}(C\com)\dto^\pi}\\
& {D^n\rto^{d\phantom{MM}}} & {Z^{n+1}(D\com)}.
\end{diagram}$$
Thus we can choose elements $h(y_i)\in C^n$ such that 
\begin{enumerate}
\item\label{hypo} $d\, h(y_i)=g(a_i)$,
\item\label{hypt} $\pi\, h(y_i)=f(y_i)$,
\end{enumerate}
for all $i\in I$. (Note that
$\pi\big(g(a_i)\big)=f\big(\iota(a_i)\big)=f(dy_i)=d\big(f(y_i)\big)$.)
By the freeness of $B^\natural$ over $A^\natural$ on $(y_i)_{i\in I}$
this defines a morphism of graded $A^\natural$-algebras
$h:B^\natural\to C^\natural$.  In particular,
$h\comp\iota=g$. Property~\ref{hypt} implies $\pi\comp h=f$.  Finally,
$hd=dh$ follows from Property~\ref{hypo}.
\end{pf}

\begin{cor}\label{qfreefib}
If $A\to B$ is a resolving morphism of differential graded algebras,
then 
\begin{equation}\label{fibr}
\shom(B,C)\longrightarrow\shom(A,C)
\end{equation}
is a fibration of simplicial sets, for all differential graded
algebras $C$.
\end{cor}
\begin{pf}
Take $X=C$ and $Y=0$ in Axiom~\ref{ax}.
\end{pf}

\begin{cor}\label{res.fib.cof}
Every resolving algebra if fibrant-cofibrant in $\AA$.
\end{cor}

Let $A\to B$ be a resolving morphism of differential graded algebras
and $A\to C$ a fixed morphism of differential graded algebras. Then we
denote the fiber of the fibration~(\ref{fibr}) over the point $A\to
C$ of $\shom(A,C)$ by $\shom_A(B,C)$.  This fiber is a
fibrant simplicial set, and we may think of it as the space of
$A$-algebra morphisms from $B$ to $C$.  

If $A\to B'\to B$ are two resolving morphisms and $A\to C$ any
fixed
morphism of differential graded algebras, then we have a fibration
$$\shom_A(B,C)\longrightarrow\shom_A(B',C)\,.$$

\subsection{Homotopies}\label{sec.homo}

The following expresses another compatibility between the closed model
category structure and the simplicial category structure on
differential graded algebras.

\begin{lem}\label{pathobject}
Let $B$ be a differential graded algebra. The canonical commutative
diagram 
\begin{equation}\label{path.o}
\vcenter{\xymatrix{
& & {B\otimes\Omega_1}\dto\\
{B\urrto\rrto_{\text{\tiny\rm diagonal}}}&&{B\times B}}}
\end{equation}
which is induced by the commutative diagram of algebraic simplices
$$\begin{diagram}
& & {\Delta^1}\dllto\\
{\Delta^0}&&{\Delta^0\amalg\Delta^0\llto\uto_{\del_0\amalg\del_1}}
\end{diagram}$$
is a path object for $B$.
\end{lem}
\begin{pf}
By the definition of path object (see \cite{simhomthe}, Section~II.1.)
we need only check that 

(i) $B\to B\otimes \Omega_1$ is a weak equivalence,

(ii) $B\otimes\Omega_1\to B\times B$ is a fibration.

Both claims reduce immediately to the case $B=k$ (by flatness of $B$
over $k$).  Then (i) is the algebraic de~Rham theorem and (ii) is
obvious.
\end{pf}

\begin{cor}
Assume that $A$ is a cofibrant differential graded algebra and $B$ a
fibrant differential graded algebra.  Two morphisms of differential
graded algebras $f,g:A\to B$ are homotopic (with respect to the closed
model category structure on $\AA$) if and only if there exists a
1-simplex $h\in \shom_1(A,B)$ such that $\del_0 h=f$ and $\del_1
h=g$. 
\end{cor}
\begin{pf}
By Corollary~1.9 of \cite[Chapter~II]{simhomthe}, the notion of
homotopy between $f$ and $g$ is well-defined. We may use the path
object of Lemma~\ref{pathobject} to check if $f$ and $g$ are
homotopic.
\end{pf}

In particular, for morphisms between resolving algebras the two
notions of homotopic are equivalent.

\begin{cor}\label{qiseq}
Let $A\to B$ be a quasi-isomorphism of resolving algebras.  Then $A\to
B$ is a homotopy equivalence.
\end{cor}
\begin{pf}
By the Theorem of Whitehead (Theorem~1.10 in
\cite[Chapter~II]{simhomthe}) $A\to B$ is a homotopy equivalence with
respect to the closed model category structure on $\AA$. 
\end{pf}

\begin{prop}
If two morphisms of resolving algebras $f,g:A\to B$ are homotopic,
then they induce the same homomorphisms on cohomology $h^\ast(A)\to
h^\ast(B)$.
\end{prop}
\begin{pf}
Let $H:A\to B\otimes\Omega_1$ be a homotopy from $f$ to $g$. Then we
have a commutative diagram of differential graded algebras
$$\begin{diagram}
& & {B\otimes\Omega_1}\dto&\\
A\urrto^H\rrto_{f\times g}&&{B\times B}& B\,,\ulto\lto
\end{diagram}$$
whose right half is our path object~(\ref{path.o}).
Applying $h^\ast$, we get the commutative diagram
$$\begin{diagram}
& & {h^\ast(B\otimes\Omega_1)}\dto&\\
h^\ast(A)\urrto^{h^\ast(H)}\ar[rr]_-{{h^\ast(f)\times
h^\ast(g)}}&&{h^\ast(B)\times h^\ast(B)}& 
h^\ast(B)\,,\ulto_{\text{\tiny isomorphism}}\lto 
\end{diagram}$$
which implies that, indeed, $h^\ast(f)=h^\ast(g)$. 
\end{pf}

\begin{cor}
Let $f:A\to B$ be a morphism of resolving algebras. If $f$ is a
homotopy equivalence, then $f$ is a quasi-isomorphism.
\end{cor}

\newcommand{\Cone}{{\cal C}}

\subsection{Derivations and differentials}\label{difcot}

Let $B\to A$ be a morphism of differential graded $k$-algebras.

\begin{defn}
For an $A$-module $M$, a {\bf $B$-derivation }$D:A\to M$ is a
homomorphism of complexes of $k$-vector spaces vanishing on $B$ and
satisfying the Leibniz rule (see~(\ref{understood}))
$$D(ab) = Da\,b+ a\,Db,$$
for all $a,b\in A$. 
\end{defn}

\begin{lem}\label{omeg}
Assume that $A$ is resolving over $B$ on the basis $(x_i)_{i\in
I}$. Then there exists a {\em universal }$B$-derivation
$\du:A\to\Omega_{A/B}$. It may be constructed as follows:

(i) As underlying graded $A^\natural$-module take
$$\Omega^\natural_{A/B}=\bigoplus_{i\in I}A^\natural\du x_i,$$
for formal generators $\du x_i$, which have the same degrees as the
$x_i$. 

(ii) Construct the unique $B$-derivation
$\du:A^\natural\to\Omega^\natural_{A/B}$, satisfying $\du(x_i)=\du
x_i$, for all $i\in I$. 

(iii) Define the differential $d$ on $\Omega^\natural_{A/B}$ by 
$$d(a\du x_i)=da\,\du x_i + (-1)^{\bar{a}} a\,\du(dx_i).$$

\noindent Then 

(i) $\Omega_{A/B}$ is an $A$-module, i.e., 
$$d(a\omega)=da\,\omega + (-1)^{\bar{a}} a\,d\omega,$$

(ii) $\du:A\to\Omega_{A/B}$ is a derivation, i.e., $d\du=\du d$,

(iii) $\du:A\to\Omega_{A/B}$ is a universal derivation.\qed
\end{lem}

\begin{cor}
If $C\to B$ is a resolving morphism of differential graded algebras,
then for every 
differential graded algebra in non-positive degrees $A$, we have a
convergent spectral sequence
$$h^q(A)\otimes_{h^0(A)} h^p\big( \Omega_{B/C}\otimes_Bh^0(A)\big)
\Longrightarrow h^{p+q}(\Omega_{B/C}\otimes_BA)\,.$$
\end{cor}
\begin{pf}
By Lemma~\ref{omeg} we know that $\Omega_{B/C}^\natural$ is free over
$B^\natural$. Since $\Omega_{B/C}$ is automatically bounded above, we
can apply Proposition~\ref{spec.se}.
\end{pf} 

\begin{defn}\label{def.imd}
The {\bf internal module of derivations }of $A$ over $B$, notation
$\Deru_B(A,M)$, is defined to be the subcomplex of $\Homu_B(A,M)$,
consisting of all elements $D:A\to M$, vanishing on $B$, and
satisfying the graded Leibniz rule
$$D(ab)=Da\,b + (-1)^{\bar{a}\bar{D}}a\,Db,$$
for all $a,b\in A$. Thus a $B$-derivation $D:A\to M$ is a
0-cocycle in the complex $\Deru_B(A,M)$. 
\end{defn}

\begin{lem}
The internal module of derivations $\Deru_B(A,M)$ is a differential
graded $A$-module. If $A$ is resolving over $B$, we have a natural
isomorphism 
$$\Deru_B(A,M)=\Homu_A(\Omega_{A/B},M).$$
Derivations correspond to homomorphisms under this isomorphism. \qed
\end{lem}

\begin{numrmk}
Let $A\to C$ be a morphism of differential graded algebras. Then
$\Deru_B(A,C)$ is a differential graded $C$-module via the action
$(cD)(a)=c\big(D(a)\big)$.
\end{numrmk}

\begin{numex}
If $A$ is a resolving algebra, we set
$$\Theta_A:=\Deru_k(A,A)=\Homu_A(\Omega_A,A)\,.$$
This has the additional structure of a differential graded Lie
algebra. The
differential on $\Theta_A$ is given by bracket with $d:A\to
A\in\Deru_k^1(A,A)$.

If $A$ is resolving over the differential graded algebra $B$, we set
$$\Theta_{A/B}=\Deru_B(A,A)=\Homu(\Omega_{A/B},A)\,.$$
\end{numex}

Now assume given morphisms of differential graded
algebras $C\to B\to A$. Let both $B$ and $A$ be resolving over $C$.
We get a homomorphism of $B$-modules $\Omega_{B/C}\to\Omega_{A/C}$
(since a $C$-derivation on $A$ restricts to a $C$-derivation on $B$),
and hence a homomorphism of $A$-modules
$\Omega_{B/C}\otimes_BA\to\Omega_{A/C}$.

\begin{lem}
If  $B\to A$ is a quasi-isomorphism, then
$$\Omega_{B/C}\otimes_BA\to\Omega_{A/C}$$ is a quasi-isomorphism.
\end{lem}
\begin{pf}
See \cite{hinich}.
\end{pf}

Assume, in addition, that $A$ is also resolving over $B$.  Then we
get a homomorphism of $A$-modules $\Omega_{A/C}\to\Omega_{A/B}$ (since
a $B$-derivation of $A$ is also a $C$-derivation).  The sequence  of
differential graded $A$-modules
\begin{equation}\label{seqdif}
0\longrightarrow \Omega_{B/C} \otimes_BA\longrightarrow \Omega_{A/C}
\longrightarrow  \Omega_{A/B}\longrightarrow 0\,
\end{equation}
is  exact.

\subsection{The cotangent complex}

Let us define the cotangent complex of a morphism between resolving
algebras.

\begin{defn}
Let $B$ and $A$ be resolving algebras over $k$, and $f:B\to A$ a
morphism. We define the {\bf cotangent complex }$L_{A/B}$ to the
differential graded $A$-module defined as the cone over the
homomorphism $\Omega_B\otimes_BA\to\Omega_A$:
$$L_{A/B}=\Cone(\Omega_B\otimes_BA\to\Omega_A).$$
\end{defn}

If $A$ is resolving over $B$, then by~(\ref{seqdif}) there is a
canonical quasi-isomorphism 
$$L_{A/B}\longrightarrow\Omega_{A/B}\,.$$

\begin{lem}\label{cot.cx}
(i) If $B\to A$ is a quasi-isomorphism, then $L_{A/B}$ is acyclic.

(ii) Given morphisms of resolving algebras $C\to B\to A$, there exists
a natural distinguished triangle of differential graded $A$-modules
$$L_{B/C}\otimes_BA\longrightarrow L_{A/C} \longrightarrow L_{A/B}
\longrightarrow L_{B/C}\otimes_BA\,[1]\,.$$

(iii) Let $A$, $B$ and $C$ be resolving algebras and  $B\to A$ a 
resolving morphism. Consider the tensor product
$$\begin{diagram}
{A\otimes_BC}& C\lto\\
A\uto & B.\lto\uto
\end{diagram}$$
So  $A\otimes_BC$ is  also a resolving algebra.
Then $L_{C/B}\otimes_BA\to L_{A\otimes_BC/A}$ and
$L_{A/B}\otimes_BC\to L_{A\otimes_BC/C}$ are quasi-isomorphisms. \qed
\end{lem}

Occasionly, we will   use the cotangent complex of a morphism
between differential graded algebras which are not resolving.
The necessary theory is developed in~\cite{hinich}.

\begin{numnote}
If $C\to B$ is a morphism of resolving algebras, then we have a
convergent spectral sequence
$$h^q(A)\otimes_{h^0(A)} h^p\big( L_{B/C}\otimes_Bh^0(A)\big)
\Longrightarrow h^{p+q}(L_{B/C}\otimes_BA)\,,$$
for every differential graded algebra in non-positive degrees $A$.\qed
\end{numnote}

\begin{defn}
Let $C\to B$ be a morphism of resolving algebras. Then we call
$$T_{B/C}=\Homu_B(L_{B/C},B)$$
the {\bf tangent complex }of $B$ over $C$.
\end{defn}

If $C\to B$ is a resolving morphism of resolving algebras, we have a
canonical quasi-isomorphism of differential graded $B$-modules
$$\Theta_{B/C}\longrightarrow T_{B/C}\,.$$

\subsubsection{Acyclicity criteria}

In the following proposition we use the differential graded algebras
$\Lambda_n$, defined for every $n>0$ as follows:

If $n$ is even, the underlying graded $k$-algebras is given by
$\Lambda_n^\natural=k[x,\xi]$, where $\deg x=-n$ and $\deg\xi=-2n-1$,
and the differential is given by $dx=0$ and $d\xi=x^2$. 

If $n$ is odd, $\Lambda_n=k[x]$, with $\deg x=-n$ and $d(x)=0$. 

Note that 
$$h^i(\Lambda_n)=\begin{cases} k & \text{if $k=0,-n$,}\\
0&\text{otherwise.}\end{cases}$$
Moreover, $\Lambda_n$ is quasi-isomorphic to $h\upst(\Lambda_n)$. 

\begin{prop}\label{eac}
Let $B\to A$ be a morphism of resolving algebras and $r\leq 0$ an
integer. Fix a morphism $A\to k$ of differential graded algebras. The
following are equivalent: 

(i) $\Deru(A,k)\to\Deru(B,k)$ is a quasi-isomorphism,

(ii) $L_{A/B}\otimes_A k$ is acyclic,

(iii) $h^r\Deru(A,C)\to h^r\Deru(B,C)$ is an isomorphism for all
(finite) resolving algebras $C$,

(iv) $h^r\Deru(A,\Lambda_n)\to h^r\Deru(B,\Lambda_n)$ is an
isomorphism for all $n>0$. 

\noindent In (iii) and (iv) the $A$-module structure on $C$ and
$\Lambda_n$ is given via $A\to k$.

If, moreover, $B\to A$ is a resolving morphism, further equivalent
conditions are

(iii$\,'$) $h^r\Deru_B(A,C)=0$, for all
(finite) resolving algebras $C$,

(iv$\,'$) $h^r\Deru_B(A,\Lambda_n)=0$, for all $n>0$. 
\end{prop}
\begin{pf}
By definition, we have a distinguished triangle
$$\Homu_A(L_{A/B},k)\longrightarrow\Deru(A,k)
\longrightarrow\Deru(B,k) \longrightarrow \Homu_A(L_{A/B},k)[1]\,.$$
Thus (i) implies that
$\Homu_A(L_{A/B},k)=\Homu_k(L_{A/B}\otimes_Ak,k)$ 
is acyclic, which implies (ii).  Let us now assume that (ii)
holds. Then we use the distinguished triangle
$$\Homu_A(L_{A/B},C)\longrightarrow\Deru(A,C)
\longrightarrow\Deru(B,C) \longrightarrow \Homu_A(L_{A/B},C)[1]\,$$
and the fact that $\Homu_A(L_{A/B},C)=\Homu_k(L_{A/B}\otimes_Ak,C)$,
to conclude that (iii) holds. The fact that (iii) implies (iv) is
trivial, because $\Lambda_n$ is a finite resolving algebra, for all
$n$. 

Finally, assume that (iv) holds. Note that we have
\begin{align*}
h^r\Deru(A,\Lambda_n)& =h^r\big(\Deru(A,k)\otimes \Lambda_n\big)\\
&=h^r\Deru(A,k)\oplus h^{r+n}\Deru(A,k)\,.
\end{align*}
Thus we may conclude that $h^\ell\Deru(A,k)\to h^\ell\Deru(A,k)$ is an
isomorphism for all $\ell\geq r$, hence for all $\ell\geq0$. Then (i)
follows, because for $\ell<0$ we have
$h^\ell\Deru(A,k)=h^\ell\Deru(B,k)=0$. 
\end{pf}

\begin{prop}\label{mea}
Let $B\to A$ be a morphism of quasi-finite resolving algebras and
$r\leq0$ and integer. The following are equivalent:

(i) $L_{A/B}$ is acyclic,

(ii) $\Deru(A,C)\to\Deru(B,C)$ is a quasi-isomorphism, for all
morphisms of differential graded algebras $A\to C$, where $C$ is a
(finite) resolving algebra,

(iii) $h^r\Deru(A,C)\to h^r\Deru(B,C)$ is an isomorphism, for all
$A\to C$ as in~(ii).

If, moreover, $B\to A$ is a resolving morphism, then a further
equivalent conditions is

(iii') $h^r\Deru_B(A,C)=0$, for all $A\to C$ as in~(ii).
\end{prop}
\begin{pf}
The claim that (i) implies (ii) follows as in the proof of
Proposition~\ref{eac} from a distinguished triangle. Then (ii) implies
(iii) trivially.  So let us assume that (iii) holds. To prove~(i), we
may assume without loss of generality that $k$ is algebraically
closed. We may
conclude from Proposition~\ref{eac} that $L_{A/B}\otimes_Ak$ is
acyclic, for all $A\to k$. This implies that $L_{A/B}$ is acyclic,
because $A$ is quasi-finite, and so $h^\ell (L_{A/B})$ is a finitely
generated $h^0(A)$-module, for all $\ell$, and hence Nakayama's lemma
applies. 
\end{pf}

\Section{\'Etale Morphisms}\label{sec.etale}

\subsection{Augmentations}

Let $A$ be a  resolving algebra. An {\em
augmentation } of $A$ is a morphism of differential graded
$k$-algebras $A\to k$. Note that every augmentation $A\to k$ is
induced (in a unique way) from a $k$-algebra morphism $h^0(A)\to
k$. Moreover, every augmentation $A\to k$ induces a morphism of
$k$-algebras $A^0\to k$. 

The {\em augmentation ideal }$\Mm$ is the kernel of $A\to k$. It is a
differential graded ideal in $A$. The augmentation ideal of $A^0$ is
the degree zero component $\Mm^0$ and the augmentation ideal of
$h^0(A)$ is $\Mm h^0(A)=\Mm^0 h^0(A)$. 

Note that any basis $(x_i)$ of $A$ defines an
augmentation $A\to k$ by the rule $x_i\mapsto 0$. Conversely, for any
augmentation of $A$, we can find a basis $(x_i)$ for $A$,
with $x_i$ in the kernel of $A\to k$, for all $i$. Such a basis is
called {\em compatible }with the augmentation $A\to k$.  If 
$(x_i)$ is a basis  compatible with the
augmentation $A\to k$, then we get an induced isomorphism of augmented
graded algebras $k[x]\to A^\natural$.

\begin{lem}\label{compstru.1}
Let $A$ be a  resolving algebra and $A\to k$ an
augmentation with ideal
$\Mm$. Then $\Mm/\Mm^2=\Omega_A\otimes_Ak$.
\end{lem}
\begin{pf}
Consider the canonical map
\begin{align*}
\Mm& \longrightarrow\Omega_A\otimes_Ak\\*
x&\longmapsto \du x\otimes 1\,,
\end{align*}
and prove that it induces an isomorphism of complexes of $k$-vector
spaces.  For this purpose it is useful to choose a basis
$(x_i)$ for $A$, compatible with the augmentation. Then
$\Omega_A\otimes_Ak$ is free as a graded $k$-vector space on 
the basis $(\du x_i\otimes 1)$, so surjectivity is clear. For
injectivity, use that  we have
$$f(x)\equiv f(0)+\sum_ix_i\frac{\del f}{\del x_i}(0)\mod \Mm^2\,,$$
for every $a=f(x)\in A$.
\end{pf}

Given an augmented resolving algebra $A\to k$, all powers $\Mm^n$ of
the augmentation ideal are differential ideals. Hence all
$\Mm^n/\Mm^{n+1}$ are complexes of $k$-vector spaces.  The direct sum
$$\gr A=\bigoplus_{n\geq0}\Mm^n/\Mm^{n+1}$$
of these complexes has an induced multiplication, making it a
differential graded algebra.  It is free:

\begin{lem}\label{compstru.2}
Let $A$ be a resolving algebra and $\Mm$ an augmentation ideal. Then
$$\gr A=S(\Mm/\Mm^2)\,,$$
as differential graded algebras.
\end{lem}
\begin{pf}
The inclusion of the subcomplex $\Mm/\Mm^2\to \gr A$ induces the
canonical morphism of differential graded algebras
$S(\Mm/\Mm^2)\to\gr A$. One checks that for every $n$ the induced
map $S^n(\Mm/\Mm^2)\to \Mm^n/\Mm^{n+1}$ is an isomorphism by
considering a basis $(x_i)$ for $A$.
\end{pf}

\subsection{Point-wise \'etale morphisms}

If $A\to k$ is a quasi-finite augmented resolving algebra, we let
$\widehat{A}^0$ 
be the completion of $A^0$ at its augmentation ideal $\Mm^0$ and 
$h^0(A)^\wedge$ the completion of $h^0(A)$ at its  augmentation ideal
$\Mm h^0(A)$. Note that we have
$h^0(A)^\wedge=h^0(A)\otimes_{A^0}\widehat{A}^0$.

We let $\widehat{A}^r$ be the completion of the $A^0$-module $A^r$ at
$\Mm^0$ and $h^r(A)^\wedge$ the completion of the $h^0(A)$-module
$h^r(A)$ at $\Mm h^0(A)$. Again, we have
$h^r(A)^\wedge=h^r(A)\otimes_{A^0}\widehat{A}^0$. The direct sum
$$\widehat{A}^\ast=\bigoplus_{r\leq0}
\widehat{A}^r=A\otimes_{A^0}\widehat{A}^0$$ 
is a differential graded algebra. We have 
$$h^r(\widehat{A}^\ast)=h^r(A\otimes_{A^0}\widehat{A}^0)=
h^r(A)^\wedge\,,$$ 
for every $r\leq0$.

If $A\to B$ is a morphism of resolving algebras, then any augmentation
$B\to k$ induces an augmentation $A\to k$. If $A$ is endowed with the
induced augmentation, then $A\to B$ is a {\em morphism }of augmented
resolving algebras.  A morphism of quasi-finite augmented resolving
algebras $A\to 
B$ induces morphisms of $k$-algebras $\widehat{A}^0\to\widehat{B}^0$
and $h^0(A)^\wedge\to h^0(B)^\wedge$. For every $r$, it induces a
homomorphism of $\widehat{A}^0$-modules
$\widehat{A}^r\to\widehat{B}^r$ and a homomorphism of
$h^0(A)^\wedge$-modules $h^r(A)^\wedge\to h^r(B)^\wedge$. Finally, it
induces a morphism of differential graded algebras
$\widehat{A}^\ast\to\widehat{B}^\ast$. 

\begin{prop}\label{pointetale}
Let $A\to B$ be a morphism of quasi-finite augmented resolving
algebras.
The following are equivalent:

(i) $L_{B/A}\otimes_Bk$ is acyclic,

(ii) $L_{B/A}\otimes_{B^0}\widehat{B}^0$ is acyclic,

(iii) $A/\Mm_A^n\to B/\Mm_B^n$ is a quasi-isomorphism for all $n$,

(iv) $\widehat{A}^\ast\to\widehat{B}^\ast$ is a quasi-isomorphism,

(v) For all $r\leq0$, the homomorphism $h^r(A)^\wedge\to
h^r(B)^\wedge$ of $h^0(A)^\wedge$-modules is bijective,

(vi) $h^0(A)\to h^0(B)$ is \'etale at $h^0(B)\to k$ and
$h^r(A)\otimes_{h^0(A)}h^0(B)\to h^r(B)$ is an isomorphism in a
Zariski neighborhood of $h^0(B)\to k$, for all $r\leq0$. 
\end{prop}
For more equivalent statements, see Proposition~\ref{eac}.
\begin{pf}
The equivalence of (iv) and (v) follows immediately from the
preceding remarks. Let us prove the equivalence of (v) and (vi):
By properties of \'etaleness for usual finite type $k$-algebras, we
know that $h^0(A)\to h^0(B)$ is \'etale at $h^0(B)\to k$ if and only
if $h^0(A)^\wedge\to h^0(B)^\wedge$ is an isomorphism of $k$-algebras.
Let us assume this to be the case. Then we have
$$h^r(A)^\wedge=h^r(A)\otimes_{h^0(A)}h^0(A)^\wedge=
\big(h^r(A)\otimes_{h^0(A)}h^0(B)\big)
\otimes_{h^0(B)}h^0(B)^\wedge$$
and
$$h^r(B)^\wedge=h^r(B)\otimes_{h^0(B)}h^0(B)^\wedge\,,$$
because $A$ and $B$ are quasi-finite. Now using the fact that for
finitely generated $h^0(B)$-modules $M$ and $N$, a homomorphism $M\to
N$ is an isomorphism in a Zariski-open neighbourhood of $h^0(B)\to k$
if and only if $(M\to N)\otimes_{h^0(B)}h^0(B)^\wedge$ is an
isomorphism of $h^0(B)^\wedge$-modules, we conclude the proof
that~(v) and~(vi) are equivalent.

The fact that (ii) implies (i) is clear.

Let us now prove that (i) implies (iii). Assume that
$L_{B/A}\otimes_Bk$ is acyclic. Then $\Omega_A\otimes_Ak\to
\Omega_B\otimes_Bk$ is a quasi-isomorphism of complexes of $k$-vector
spaces. By Lemma~\ref{compstru.1}, $\Mm_A/\Mm_A^2\to\Mm_B/\Mm_B^2$ is
a quasi-isomorphism. Since taking symmetric powers preserves
quasi-isomorphisms, Lemma~\ref{compstru.2} implies that we have a
quasi-isomorphism $\Mm_A^n/\Mm_A^{n+1}\to\Mm_B^n/\Mm_B^{n+1}$, for all
$n\geq0$. By induction, this implies that $A/\Mm_A^n\to B/\Mm_B^n$ is
a quasi-isomorphism, for all $n$, proving (iii). 

Assuming (iii), we will now prove (iv).
Let us fix $r\leq0$ and consider the limit
$$\projectlim(A/\Mm^n)^r\,.$$
(The superscript $n$ denotes a power, the superscript $r$ denotes the
component of degree $r$.) Note that this limit is equal to
$\widehat{A}^r$, because the topology on $A^r$ defined by the
descending sequence of subspaces $(\Mm^n)^r$ is equal to the
$\Mm^0$-adic topology: for every $n$ we have
$$(\Mm^0)^nA^r\subset(\Mm^n)^r\subset(\Mm^0)^{n-r}A^r\,.$$
Thus Lemma~\ref{lim} implies that we
have an isomorphism $h^r(\widehat{A}^\ast)\to
h^r(\widehat{B}^\ast)$. 

Finally, let us prove that (iv) implies (ii).
\begin{align*}
\widehat{A}^\ast\to \widehat{B}^\ast\quad\text{qis}\quad & \Longrightarrow  
L_{\widehat{A}^\ast}\otimes_{\widehat{A}^\ast}\widehat{B}^\ast
\to  L_{\widehat{B}^\ast}\quad\text{qis}\\  
&\Longrightarrow 
\Omega_{A}\otimes_{A}\widehat{B}^\ast\to
\Omega_{B}\otimes_B\widehat{B}^\ast 
\quad\text{qis}\,.
\end{align*} 
For the necessary facts about cotangent complexes we refer to
\cite{hinich}.  
\end{pf}

\begin{lem}\label{lim}
Let $M$ be a complex of $k$-vector spaces and $N_n\subset
M$ a descending sequence of subcomplexes. Then for every $r$
there is a natural exact sequence
\begin{multline*}
\prod_n h^{r-1}(M/N_n)\stackrel{\del}{\longrightarrow}
\prod_n h^{r-1}(M/N_n)\stackrel{\sigma}{\longrightarrow}\\
\stackrel{\sigma}{\longrightarrow}
h^r(\projectlim M/N_n)\stackrel{\iota}{\longrightarrow}
\prod_n h^{r}(M/N_n)\stackrel{\del}{\longrightarrow}
\prod_n h^{r}(M/N_n)\,.
\end{multline*}
Here $\projectlim M/N_n$ denotes the componentwise projective limit
of complexes.
\end{lem}
\begin{pf}
This is straightforward to check from 
the definitions:
The map $\del$ maps the sequence $a=(a_n)$ to $\del a$, with $(\del
a)_n= a_n-a_{n+1}$. The map $\sigma$ sends the sequence $a$ to
$\sigma(a)$ with $\sigma(a)_n=\sum_{i=1}^{n-1} da_i$.  The map $\iota$
sends the sequence $a$ to $a$.  One can also start with the short
exact sequence of 
complexes whose degree $r$ term is
$$\ses{\projectlim (M/N_n)^r}{\iota}{\prod_n (M/N_n)^r}{\del}{\prod_n
(M/N_n)^r}\,,$$
and then pass to the long exact cohomology sequence.
\end{pf}

\begin{defn}\label{defpointetale}
Let $A\to B$ be a morphism of quasi-finite resolving algebras. Let
$B\to k$ be an augmentation.
If the equivalent conditions of Corollary~\ref{pointetale} are
satisfied, then we call $A\to B$ {\bf \'etale }at the augmentation
$B\to k$. 
\end{defn}

If $k$ is not algebraically closed, there might not exist very many
augmentations.  Thus the following generalization:  Let $K$ be a
finite extension field of $k$.  A morphism of differential graded
algebras $A\to K$ is called a {\em $K$-valued }augmentation of $A$. 

\begin{cor}
Let $A\to B$ be a morphism of quasi-finite resolving algebras and
$B\to K$ a $K$-valued augmentation of $B$. Then the following are
equivalent:

(i) $L_{B/A}\otimes_B K$ is acyclic,

(ii) $h^0(A)\to h^0(B)$ is \'etale at $h^0(B)\to K$ and
$h^r(A)\otimes_{h^0(A)}h^0(B)\to h^r(B)$ is an isomorphism in a
Zariski neighborhood of $h^0(B)\to K$, for all $r\geq0$. 

(iii) $A\otimes_kK\to B\otimes_kK$ is \'etale at the augmentation
$B\otimes_k K\to K$.
\end{cor}

\subsection{\'Etale morphisms}

\begin{cor}\label{condetale}
Let $A\to B$ be a morphism of quasi-finite resolving algebras.
The following are equivalent:

(i) $L_{B/A}$ is acyclic,

(ii) $h^0(A)\to h^0(B)$ is \'etale and
$h\upst(A)\otimes_{h^0(A)}h^0(B)\to h\upst(B)$ is an isomorphism,

(iii) $A\to B$ is \'etale at every $K$-valued augmentation $B\to K$,
for all finite extension $K/k$.
\end{cor}
For more equivalent statements, see Proposition~\ref{mea}.
\begin{pf}
Proposition~\ref{pointetale} implies directly that 
Conditions (ii) and (iii) are equivalent.  To deal with Condition~(i),
we may assume that $k$ is algebraically closed. We notice 
that $h^r(L_{B/A})$ is a finitely generated $h^0(B)$-module, for every
$r$. Thus,  $L_{B/A}$ is acyclic, if and only if
$L_{B/A}\otimes_{B^0}\widehat{B}^0$ is acyclic for all augmentations
$h^0(B)\to k$.
\end{pf}

\begin{defn}\label{defetale}
If the equivalent conditions of Corollary~\ref{condetale} are
satisfied, we call the morphism of quasi-finite resolving algebras
$A\to B$ {\bf \'etale}. 
\end{defn}

\begin{cor}\label{qiscondition}
Let $A\to B$ be a morphism of quasi-finite resolving algebras. Then
$A\to B$ is a 
quasi-isomorphism if and only if $h^0(A)\to h^0(B)$ is an isomorphism
and $L_{B/A}$ is acyclic.\qed
\end{cor}

\begin{prop} \label{propetale}
Let $A\to B\to C$ be morphisms of quasi-finite resolving algebras,
where $A\to B$ is \'etale. Then $B \to C$ is \'etale if and only if
$A\to B$ is \'etale.\qed
\end{prop}

\begin{prop}
Let $A\to A'$ be an \'etale morphism of quasi-finite resolving
algebras.  Let $M$ be a differential graded $A$-module, bounded from
above and such that $M^\natural$ is flat over $A^\natural$.  Then we
have $$h^p(M\otimes_AA')=h^p(M)\otimes_{h^0(A)}h^0(A')\,,$$
for all $p$. 
\end{prop}
\begin{pf}
Since usual \'etale morphisms of finite type $k$-algebras are flat,
$h^0(A)\to h^0(A')$ is flat.  Since flatness is preserved under base
change, this implies that $h\upst(A)\to h\upst(A')$ is flat. Hence 
$$\Tor^{h\upst(A)}_i\big(h\upst(M),h\upst(A')\big)=0\,,$$
for all $i>0$.  By the Eilenberg-Moore spectral sequence (see
\cite{DerQuot}) this implies that
$$h\upst(M\otimes_AA') =h\upst(M)\otimes_{h\upst(A)}h\upst(A')
=h\upst(M)\otimes_{h^0(A)}h^0(A')\,,$$  
which is what we wanted to prove.
\end{pf}

The following special case of this proposition will be essential for
descent theory. Once we have defined {\em perfect resolving morphism
}(see Section~\ref{sec.perfect}), it is clear that this corollary
generalizes to the case that $C\to B$ is a perfect resolving morphism.

\begin{cor}\label{prep.desc}
Let $C\to B\to A \to A'$ be morphisms of quasi-finite resolving
algebras over $k$.  Assume that $C\to B$ is a finite resolving
morphism and $A\to A'$ is 
\'etale.  Then we have
$$h^p\big(\Deru_C(B,A')\big)
=h^p\big(\Deru_C(B,A)\big)\otimes_{h^0(A)}h^0(A')\,,$$ 
for all $p$.\qed
\end{cor}

\begin{defn}\label{def.op.imm}
We call $A\to B$ an {\bf open immersion}, if 

(i) $A\to B$ is \'etale,

(ii) $\spec h^0(B)\to\spec h^0(A)$ is an open immersion of affine
$k$-schemes.
\end{defn}

\subsection{Completions}

Here we would like to characterize \'etaleness in terms of
completions.  This section is not used in the rest of the paper.  We
include it here, even
though it contains a forward reference to the notion of perfect
resolving algebra.

Let $A$ be a quasi-finite  resolving algebra and $A\to k$ be an
augmentation.  Let 
$\Mm$ be the augmentation ideal. 
the various quotients 
$A/\Mm^{n}$,  form an inverse
system of  
differential graded $k$-algebras.  We define
$$\widehat{A}=\projectlim A/\Mm^n,$$
as a projective limit of $k$-algebras and call it the {\em completion
}of $A$ at $\Mm$ or at the augmentation $A\to k$.   The completion
$\widehat{A}$ inherits the structure of differential $k$-algebra, but it
will lose its grading. 

The differential algebra $\widehat{A}$ comes with a natural
filtration.  The components of the filtration are the differential
ideals $\widehat{\Mm}_n=\ker(\widehat{A}\to A/\Mm^n)$ and the
associated graded pieces are the $\widehat{A}$-modules
$\widehat{\Mm}_n/\widehat{\Mm}_{n+1}$. Note that
$$\widehat{\Mm}_n/\widehat{\Mm}_{n+1}=\Mm^n/\Mm^{n+1}\,.$$

Choose a coordinate system $(x_i)$ for $A$, compatible with the
augmentation. Let $\langle x\rangle$ be the graded vector space
generated by the $x_i$ and let $k[[x]]=\widehat{S}\langle
x\rangle=\prod_{n\geq0}S^n\langle x\rangle$ be
the completion of the symmetric algebra at the augmentation ideal
$(x)\subset k[x]$. We may view the elements of
$k[[x]]=\widehat{A}^\natural$ as formal power series in the variables
$x_i$.

Note that
$$\widehat{A}^r=\{\hat{a}\in\widehat{A}\st\text{for all $n$, the
projection of $\hat{a}$ into $A/\Mm^n$ has degree $r$}\}\,.$$
Thus we have an inclusion of differential algebras
$$\widehat{A}^\ast\subset \widehat{A}\,,$$
which is not an isomorphism unless $A=k$.

From the various projections
$\widehat{A}\to\widehat{A}^r$, we get a canonical inclusion
\begin{equation}\label{incl}
\xymatrix@C=15pt{
{{\phantom{\displaystyle\prod_{r\leq0}}}\widehat{A}}\,\,
\ar@<1ex>@{^{(}->}[r]& 
\,{\displaystyle\prod_{r\leq0}\widehat{A}^r}\,.}
\end{equation}
This is to be understood as a morphism of differential
$k$-algebras. Multiplication on the right hand side is given by
considering elements as formal infinite sums
$\hat{a}=\sum_{r\leq0}\hat{a}^r$
and multiplying using the distributive law. The differential on the
right had side is given by
$d\hat{a}=\sum d(\hat{a}^r)$.
Note that
$d(\hat{a}^r)=(d\hat{a})^{r+1}$, for all $\hat{a}\in\widehat{A}$.

The image of (\ref{incl}) consists of all formal series
$\hat{a}=\sum_{r\leq0}\hat{a}^r$ which converge in the $\Mm$-adic
topology, i.e., such that 
$$\forall n\quad\exists R\quad\forall r>R\colon\quad
\hat{a}^r\in\widehat{\Mm}_n\,.$$
This condition is automatically satisfied if $A$ a finite resolving
algebra. So in this case, (\ref{incl}) is an isomorphism of
differential algebras.

If we take cohomology of~(\ref{incl}), we obtain an algebra morphism
\begin{equation}\label{incl.f}
h(\widehat{A})\longrightarrow
\prod_{r\leq0}h^r(\widehat{A}^\ast)\,,
\end{equation}
which is an isomorphism in the finite case.

\begin{them}\label{et.com.crit}
Let $A\to B$ be a morphism of quasi-finite resolving algebras and
$B\to k$ an augmentation.  If $A\to B$ is \'etale at $B\to k$, then we
have a quasi-isomorphism of completions $\widehat{A}\to
\widehat{B}$. The converse is true if $A$ and $B$ are perfect.
\end{them}
\begin{pf}
Assume that $A\to B$ is \'etale at $B\to k$. Then for all $n$ we have
that $A/\Mm^n\to B/\Mm^n$ is a quasi-isomorphism. 
We now employ the non-graded version of Lemma~\ref{lim}: Let $M$ be a
differential $k$-vector space and $N_n\subset M$ a descending sequence
of subspaces respecting the differential. Then we have an exact
triangle of $k$-vector spaces
$$\xymatrix@C=0pt{
{\displaystyle\prod_n h(M/N_n)}\ar@<1ex>[rr]^{\del} &&
{\displaystyle\prod_n 
h(M/N_n)}\dlto^{\sigma}\\ 
& {h(\projectlim M/N_n)}\ulto^{\iota} &}$$
where the maps are defined by the same formulas as those of
Lemma~\ref{lim}. We may obtain this exact triangle from the short
exact sequence of differential vector spaces
$$\ses{\projectlim M/N_n}{}{\prod M/N_n}{\del}{\prod M/N_n}\,.$$

Applying this to our situation, we obtain the desired result that
$h(\widehat{A})\to h(\widehat{B})$ is an isomorphism.

Conversely, assume that we have a quasi-isomorphism $\widehat{A}\to
\widehat{B}$. To prove that $A\to B$ is \'etale at $B\to k$, we may
localize $A$ and $B$ and thus assume that they are
finite. Then
we use~(\ref{incl.f}) to 
conclude that we have an isomorphism
$$\prod_{r\leq0} h^r(\widehat{A}^\ast)\longrightarrow\prod_{r\leq0}
h^r(\widehat{B}^\ast)\,,$$
which implies that $h^r(\widehat{A}^\ast)\to h^r(\widehat{B}^\ast)$ is
bijective for all $r$, hence that $\widehat{A}^\ast\to
\widehat{B}^\ast$ is a quasi-isomorphism.
\end{pf}

\begin{rmk}
Note that unlike the non-differential graded case, we do not have that
$\widehat{A}=\widehat{S}(\Mm/\Mm^2)$. By assumption on $A$, we can
always find a subcomplex $V\subset \Mm$, such that $V\to \Mm/\Mm^2$ is
an isomorphism and $A^\natural=S(V)^\natural$. Hence we also have that
the filtered algebras $\widehat{A}^\natural$ and
$\widehat{S}(\Mm/\Mm^2)^\natural$ are isomorphic. But the
differentials are different, as soon as the differential $d$ on $A$
has any higher order~($\geq2$) terms.

We may say that the tangent cone and the tangent space are isomorphic,
but the completion of the differential graded algebra itself is not
isomorphic to the tangent space.  Thus quasi-free differential graded
algebras have some, but not other properties of smooth
non-differential graded algebras.
\end{rmk}

\subsection{Local structure of \'etale morphisms}\label{sec.loc.str}

Given a differential graded algebra $A$, assume that the differential
graded algebra $B$ is quasi-free over $A$ on the basis
$x_1,\ldots,x_r$ in degree $0$ and $\xi_1,\ldots,\xi_s$ in degree
$-1$. Denote $d\xi_i$ by $f_i\in A^0[x_1,\ldots,x_r]$.  In this case
we write
$$B=A[x_1,\ldots,x_r]\{\xi_1,\ldots,\xi_s\}/d\xi_i=f_i.$$

\begin{lem}\label{standard}
Let $A$ be a quasi-finite resolving algebra. Let
$f_1,\ldots,f_r\in A^0[x_1,\ldots,x_r]$ be polynomials such that
$\det({\del f_i\over \del x_j})$ is a unit in
$h^0(A)[x_1,\ldots,x_r]/(f_1,\dots,f_r)$.  Then 
$$A\longrightarrow
A[x_1,\ldots,x_r]\{\xi_1,\ldots,\xi_r\}/d\xi_i=f_i$$ is
\'etale.
\end{lem}
\begin{pf}
Let $B=A[x]\{\xi\}/dx=f$.
By usual facts about \'etale morphisms between $k$-algebras, we know
that $h^0(A)\to h^0(B)$ is \'etale.  The assumption on the Jacobian
immediately implies that $L_{B/A}\otimes_B K$ is acyclic, for all
$K$-valued augmentations $B\to K$.
\end{pf}

\begin{defn}\label{st.et}
We call an \'etale morphism
$$A\longrightarrow
A[x_1,\ldots,x_r]\{\xi_1,\ldots,\xi_r\}/d\xi_i=f_i$$
as in Lemma~\ref{standard} a {\bf standard \'etale }morphism.
\end{defn}

\begin{note}
A composition of standard \'etale morphisms is standard
\'etale.
\end{note}

Let $g\in A^0$ and consider $f(x)=xg-1\in A^0[x]$.  Then ${\del f\over
\del x}=g$ is a unit in $A^0[x]/xg-1=A^0_g$.  Thus $A\to
A[x]\{\xi\}/d\xi=xg-1$ is \'etale.  In fact, it is an open immersion,
as it induces $h^0(A)\to h^0(A)_g$ on the $h^0$-level.  We will
abbreviate $A[x]\{\xi\}/d\xi=xg-1$ by $A_{\{g\}}$.

\begin{defn}
An open immersion $A\to A_{\{g\}}$ is called an {\bf elementary open
immersion}.
\end{defn}

\begin{prop}\label{loc.str}
Let $A\to B$ be a morphism of quasi-finite resolving algebras which is
\'etale in a Zariski neighborhood of the augmentation $h^0(B)\to k$.
Then there exists a commutative diagram of differential graded
algebras
$$\xymatrix{
A\rto\dto_{\text{\rm standard \'etale}} & B\dto^{\text{\rm elem.\
open}}\\
B''\rto^{\text{\rm qis}} & B'}$$
and an augmentation $B'\to k$ compatible with the given augmentation
$B\to k$.
\end{prop}

\begin{pf}
We use the local structure theory of usual \'etale morphisms: there exists
$g\in h^0(B)$, further polynomials $f_1,\ldots,f_p\in
h^0(A)[x_1,\ldots,x_p]$ and an $h^0(A)$-isomorphism
$h^0(A)[x]/(f)\to h^0(B)_g$, such that $\det(\frac{\del f}{\del x})$
is a unit in $h^0(A)[x]/(f)$. Moreover, we may assume that $g$ does
not map to zero under the augmentation $h^0(B)\to k$ and that $A\to B$
is \'etale over $h^0(B)_g$

We choose a lifting $g\in B^0$ and consider the
localization $B\to B_{\{g\}}=B[y]\{\eta\}/(d\eta=1-yg)$. 
We define an augmentation $B_{\{g\}}\to k$ by sending $y$ to the
inverse of the image of $g$ in $k$ and $\eta$ to
$0$ and making sure that it restricts to the given augmentation on
$B$. 
Note that $A\to B_{\{h\}}$ is \'etale.

Now lift $f_1,\ldots,f_p$ to elements of $A^0[x_1,\ldots,x_p]$,
denoted by the same letters. Also lift the  images of $x_1,\ldots,x_p$
in $h^0(B_{\{g\}})$ to elements $b_1,\ldots,b_p\in
B_{\{g\}}^0$. Because $f_i(b)$ represents $0$ in $h^0(B_{\{g\}})$ we
may also choose $\beta_i\in B_{\{g\}}^{-1}$ such that
$d\beta_i=f_i(b)$, for all $i=1,\ldots,p$. 

Having made these choices, we can define a morphism of differential
graded algebras 
\begin{align*}
A[x]\{\xi\}/(d\xi=f) & \longrightarrow B_{\{g\}}\\*
x_i&\longmapsto b_i\\*
\xi_i&\longmapsto \beta_i\,.
\end{align*}
Note that this is a quasi-isomorphism by
Corollary~\ref{qiscondition}. Thus  the commutative diagram 
$$\xymatrix{
A\dto\rto & B\dto\\
A[x]\{\xi\}/(d\xi=f)\rto^-{\text{qis}}& B_{\{g\}}}$$
finishes the proof.
\end{pf}

\begin{cor}\label{mailem}
Let $A\to B$ be an \'etale morphism of quasi-finite resolving
algebras. Then there exists 
an integer $n$, and for every $i=1,\ldots,n$ a commutative
diagram of differential graded algebras
$$\begin{diagram}
A\dto_{\text{\rm standard \'etale}} \rto & B\dto^{\text{\rm elem.\
open}}\\ 
B_i'\rto_{\text{\rm qis}}
& B_i
\end{diagram}$$
such that 

(i) all $B\to B_i$ are elementary open immersions and $\coprod_i\spec
h^0(B_i)\to\spec h^0(B)$ is surjective,

(ii) all $A\to B'_i$ are standard \'etale,

(iii) all $B'_i\to B_i$ are quasi-isomorphisms.
\end{cor}
\begin{pf}
Immediate from Proposition~\ref{loc.str}.
\end{pf}

\Section{Perfect resolving algebras}\label{sec.perfect}

Perfect resolving algebras are between finite and quasi-finite
resolving algebras.  Their importance lies in the fact that they
provide us with a good notion of affine differential graded schemes. 

\begin{defn}\label{def.per}
A {\bf perfect }resolving algebra is a quasi-finite resolving algebra
$A$, such that $\Omega_A\otimes_A h^0(A)$ is a perfect complex of
$h^0(A)$ modules. 

If $\Omega_A\otimes_A h^0(A)$ has  perfect amplitude contained in
$[-N,0]$, then we say that $A$ has {\bf  perfect amplitude }$N$. 

The full sub-2-category of $\RR$, consisting of perfect resolving
algebras will be denoted by $\RR_\perf$.
\end{defn}

Recall that a complex $M$ of $h^0(A)$-modules is {\em perfect }if
Zariski locally in $h^0(A)$, there exists a finite complex of finite
rank free modules $E$, and a quasi-isomorphism $E\to M$. If all
locally defined $E$ can be chosen such that $E^i=0$ for
$i\not\in[a,b]$, then $M$ is of {\em perfect amplitude }contained in
$[a,b]$. For details on perfect complexes see Expos\'es~I and~II
of~\cite{sga6}. 

For example, any finite resolving algebra $A$ is perfect, because if
$x_1,\ldots,x_n$ is a basis for $A$, then $\du x_1\otimes
1,\ldots,\du x_n\otimes 1$ is a basis for $\Omega_A\otimes_A
h^0(A)$.  If $A$ has amplitude $N$ as a finite resolving algebra, it
has perfect amplitude $N$. The converse is true `locally', if
$N\geq2$; see Theorem~\ref{loc.fin}.

\begin{defn}
A {\bf perfect }resolving morphism of differential graded algebras
$C\to B$ is a quasi-finite resolving morphism such that
$\Omega_{B/C}\otimes_Bh^0(B)$ is a perfect complex of
$h^0(B)$-modules.  

If $\Omega_{B/C}\otimes_B h^0(B)$ has  perfect amplitude contained in
$[-N,0]$, then $C\to B$ has {\bf  perfect amplitude }$N$. 
\end{defn}

\begin{defn}\label{depeqf}
A morphism $C\to B$ of quasi-finite resolving algebras is called {\bf
perfect}, if $L_{B/C}\otimes_Bh^0(B)$ is a perfect complex of
$h^0(B)$-modules. 

If $L_{B/C}\otimes_Bh^0(B)$ has perfect amplitude contained in
$[-N,0]$, we call $C\to B$ of {\bf perfect amplitude }$N$.
\end{defn}

For example, any \'etale morphism of quasi-finite resolving algebras
if perfect, of perfect amplitude $0$. Any morphism between finite
resolving algebras is perfect. If $A$ has amplitude $N$ and $B$ has
amplitude $M$, then $A\to B$ has perfect amplitude $\max(N+1,M)$. 
More precisely, we have:

\begin{numrmk}\label{formnp}
Let $C\to B\to A$ be morphisms of differential graded algebras.
Assume either that both $C\to B$ and $B\to A$ are quasi-finite
resolving morphisms, or that all three of $C$, $B$, $A$ are
quasi-finite resolving algebras. 
$$\xymatrix{C\rto^f\drto_h& B\dto^g\\& A}$$
If any two of $f$, $g$, $h$ are perfect, then so is the third. If we
denote the amplitudes of $f$, $g$ and $h$, by $M$, $N$ and $P$,
respectively, then we have
\begin{align*}
P&=\max(M,N)\\
N&=\max(M+1,P)\\
M&=\max(N-1,P)
\end{align*}
For a reference, see Compl\'ement~4.11 in Expos\'e~I of~\cite{sga6}.
\end{numrmk}

The first basic result about perfect morphisms is that we have a
convergent spectral sequence for the tangent complex. 

\begin{prop}\label{ss.th}
Let $C\to B$ be a perfect resolving morphism of differential graded
algebras and $B\to A$ a morphism to a differential graded algebra $A$.
Assume that $A$ is concentrated in non-positive degrees (it is
sufficient that the cohomology of $A$ be concentrated in non-positive
degrees).   Then, for all $p$ sufficiently
large, we have that $h^p\big(\Theta_{B/C}\otimes h^0(A)\big)=0$. Thus
there 
is a
convergent third third quadrant spectral sequence
$$E_2^{p,q}=h^q(A)\otimes_{h^0(A)} h^p\big(\Theta_{B/C}\otimes_B
h^0(A)\big)\Longrightarrow h^{p+q}(\Theta_{B/C}\otimes_BA)\,.$$
Rewriting in terms of derivations yields 
\begin{equation}\label{ss.th.t}
E_2^{p,q}=h^q(A)\otimes_{h^0(A)}
h^p\Deru_C\big(B,h^0(A)\big)\Longrightarrow
h^{p+q}\Deru_C(B,A)\,.
\end{equation}
More precisely, if $C\to B$ has perfect amplitude $N$, then
$\Theta_{B/C}\otimes_B h^0(A)$ has perfect amplitude contained in
$[0,N]$, and hence 
$h^p\big(\Theta_{B/C}\otimes_B h^0(A)\big)=0$, for all $p>N$. 

Let $C\to B$ be a perfect morphism of quasi-finite resolving
algebras. Then we have a convergent third quadrant spectral sequence
$$E_2^{p,q}=h^q(A)\otimes_{h^0(A)} h^p\big(T_{B/C}\otimes_B
h^0(A)\big)\Longrightarrow h^{p+q}(T_{B/C}\otimes_BA)\,.$$
If $C\to B$ has perfect amplitude $N$, then $T_{B/C}\otimes_B
h^0(A)$ has perfect amplitude contained in $[0,N]$ and
$h^p\big(T_{B/C}\otimes_B 
h^0(A)\big)=0$, for all $p>N$. 
\end{prop}
\begin{pf}
We note that 
$$\Theta_{B/C}\otimes_Bh^0(A)=
\Homu_{h^0(A)}\big(\Omega_{B/C}\otimes_Bh^0(A),h^0(A)\big)\,.$$
Thus the perfection of $\Omega_{B/C}\otimes_Bh^0(A)$ implies that of
$\Theta_{B/C}\otimes_Bh^0(A)$. Then the proposition follows.
\end{pf}

\begin{cor}
If $C\to B$ is a perfect resolving morphism and $A$ is concentrated in
non-positive degrees, then $h^\ell\Deru_C(B,A)$ is a finite rank
$h^0(A)$-module. 
\end{cor}

\subsubsection{Application to the behaviour of $\Deru$ with respect to
homotopies} 

Let $C\to B$ be a perfect resolving morphism of differential graded
algebras.
Let $f:B\to A$ and $g:B\to A$ be morphisms, where the cohomology of
$A$ is concentrated in non-positive degrees.  Let
$\theta:f\Rightarrow g$ be a homotopy, i.e., a morphism $\theta:B\to
A\otimes\Omega_1$, such that $\del_0\theta=f$ and $\del_1\theta=g$. 
$$\xymatrix{
B\rtwocell^f_g{_\theta} & A}$$
In
this situation, there are two $B$-module structures on $A$, one given
by $f$, the other by $g$. Let us denote them by ${_fA}$ and ${_gA}$.
The algebra $A\otimes\Omega_1$ has three $B$-module structures.  Let
us denote the one induced by $\theta$ by ${_\theta
A\otimes\Omega_1}$. 

We get induced homomorphisms of $B$-modules
$$\xymatrix@R=.5pc{
&& {\Deru(B,{_fA})}\\
{\Deru(B,{_\theta A\otimes\Omega_1})}\urrto^{{\del_0}\lst}
\drrto_{{\del_1}\lst} \\
&& {\Deru(B,{_gA})}}$$
By Proposition~\ref{ss.th}, both of these are quasi-isomorphisms.
Thus we can make the following definition:
 
\begin{defn}\label{can-iso}
The isomorphism of $h^0(B)$-modules
$$h_\ell \Deru(B,{_fA})\stackrel{\theta\lst}{\longrightarrow}
h_\ell\Deru(B,{_gA})$$ 
obtained as the composition of $({\del_0}\lst)^{-1}$ with
${\del_1}\lst$, is called the {\bf canonical }isomorphism induced by
$\theta$, and is denoted by $\theta\lst$.
\end{defn}

Note that the canonical isomorphism induced by $\theta$ depends only
on the homotopy class of $\theta$ and is
functorial: $\theta\lst\eta\lst=(\theta\eta)\lst$. 

\begin{rmk}
In fact,  $\theta\lst$ is entirely independent of
$\theta$.
\end{rmk}

\subsection{Local finiteness}

We will prove two fundamental results on perfect resolving algebras.
The first one says that every perfect resolving algebra is locally
finite:

\begin{them}\label{loc.fin}
Let $C\to B$ be a perfect  resolving morphism of
quasi-finite resolving 
algebras. For every $K$-valued augmentation $B\to K$, there exists a
$g\in B^0$, such that $g(K)\not=0$ and a resolution $C\to A\to
B_{\{g\}}$, of $C\to B_{\{g\}}$, where $C\to A$ is a finite resolving
morphism. 
$$\xymatrix{
C\rto\dto_{\text{\rm finite}} & B\dto\\
A\rto^-{\text{\rm qis}} & B_{\{g\}}}$$
If $C\to B$ has perfect amplitude $N$, then we can choose
$A$ such that $C\to A$ has amplitude $\max(2,N)$. 
\end{them}
\begin{pf}
Since $C\to B$ is a resolving morphism, we have a quasi-isomorphism
$L_{B/C}\to\Omega_{B/C}$, and so we know that
$\Omega_{B/C}\otimes_Bh^0(B)$ is perfect.  Suppose that $n\geq2$ is an
integer such that $\Omega_{B/C}\otimes_Bh^0(B)$ has perfect amplitude
contained in $[-n,0]$. Then
$\tau_{\geq-n}\big(\Omega_{B/C}\otimes_Bh^0(B)\big)$ has  perfect
amplitude contained in $[-n,0]$, too, because
$\Omega_{B/C}\otimes_Bh^0(B)\to
\tau_{\geq-n}\big(\Omega_{B/C}\otimes_Bh^0(B)\big)$ is a
quasi-isomorphism.
Now, since $\tau_{\geq-n}\big(\Omega_{B/C}\otimes_Bh^0(B)\big)^i$ is a
free $h^0(B)$-module, for all $i>-n$ (a basis if provided by $(\du
x_j\otimes 1)$, where $(x_j)$ is the degree $i$ part of a basis for 
$B$ over $C$), it follows that
\begin{multline*}
\tau_{\geq-n}\big(\Omega_{B/C}\otimes_Bh^0(B)\big)^{-n}=\\
=\cok\bigg(
\big(\Omega_{B/C}\otimes_Bh^0(B)\big)^{-n-1} \longrightarrow
\big(\Omega_{B/C}\otimes_Bh^0(B)\big)^{-n}\bigg)
\end{multline*}
is a locally free $h^0(B)$-module (for a reference, see Lemme~4.16,
Expos\'e~I in~\cite{sga6}).  Let $h^0(B)_g$ be a Zariski neighborhood
of the $K$-valued point $h^0(B)\to K$ of $\spec h^0(B)$ over which the
above cokernel is free. Lifting $g$ to $B^0$, we get the elementary
open immersion $B\to B_{\{g\}}$ through which $B\to K$ factors. Then
$\Omega_{B_{\{g\}}/C}\otimes_{B_{\{g\}}}h^0(B_{\{g\}})$ has perfect
amplitude contained in $[-n,0]$, because $n\geq1$, and 
$$\tau_{\geq-n}\big(\Omega_{B_{\{g\}}/C}
\otimes_{B_{\{g\}}}h^0(B_{\{g\}})\big)^{-n}=
\tau_{\geq-n}\big(\Omega_{B/C}\otimes_Bh^0(B)\big)^{-n}
\otimes_{h^0(B)}h^0(B)_g$$
is a free $h^0(B_{\{g\}})=h^0(B)_g$-module.  Thus, to simplify
notation, we may replace $B$ by $B_{\{g\}}$, and assume that
$\tau_{\geq-n}\big(\Omega_{B/C}\otimes_Bh^0(B)\big)^{-n}$ is free.

Choose elements $b_1,\ldots,b_r\in B^{-n}$, such that the images of
$\du b_j\otimes 1\in\big(\Omega_{B/C}\otimes_B h^0(B)\big)^{-n}$ in
$\tau_{\geq-n}\big(\Omega_{B/C}\otimes_Bh^0(B)\big)^{-n}$ give a
basis. Then define $A=B_{(n-1)}[\xi_1,\ldots,\xi_r]$, where
$\xi_1,\ldots,\xi_r$ are formal variables of degree $-n$ and
$B_{(n-1)}\subset B$ is the differential graded $C$-subalgebra
generated over $C$ by the part of degree $>-n$ of a basis for $B$
over $C$. Set $d\xi_j=db_j$, which is allowed, because
$db_j\in Z^{1-n} B_{(n-1)}$. Extend the inclusion $B_{(n-1)}\subset B$
to a morphism of differential graded algebras $A\to B$ by
$\xi_j\mapsto b_j$. Note that $C\to B_{(n-1)}$, and hence $C\to A$ are
finite resolving morphisms. The amplitude of $A$ is $n$. 

Now, all that is left, is to check that $A\to B$ is a
quasi-isomorphism. We use Corollary~\ref{qiscondition}: $h^0(A)\to
h^0(B)$ is an isomorphism since $n\geq2$. Moreover, we have a
commutative diagram of complexes of $h^0(B)$-modules
$$\xymatrix{
\Omega_{A/C}\otimes_Ah^0(B)\rto\drto_{\text{\rm isom.}} &
\Omega_{B/C}\otimes_Bh^0(B) \dto^{\text{\rm qis}}\\
& \tau_{\geq-n}\big(\Omega_{B/C}\otimes_Bh^0(B)\big)\,, }$$
where the diagonal map is an isomorphism by  construction. This
proves that $L_{B/A}\otimes_Bh^0(B)$, and hence $L_{B/A}$, is acyclic.
\end{pf}

\begin{cor}
Let $C\to B$ be a perfect morphism of quasi-finite resolving
algebras. Then for every $K$-valued augmentation $B\to K$ there exists
an open immersion $B\to B'$, such that $B\to K$ factors through $B\to
B'$, and a commutative diagram 
$$\xymatrix{
C\rto\dto_{\text{\rm finite}} & B\dto^{\text{\rm open imm.}}\\
A\rto^-{\text{\rm qis}} & B'}$$
where $C\to A$ is a finite resolving morphism and $A\to B'$ a
quasi-isomorphism. 
If $C\to B$ has perfect amplitude $N$, then $C\to A$ can be chosen to
have amplitude $\max(2,N)$.
\end{cor}
\begin{pf}
Start by choosing (using Proposition~\ref{res.ex}) a quasi-finite
resolution $C\to B''\to B$ of $C\to 
B$, and choose a homotopy inverse for $B''\to B$, to obtain the
homotopy commutative 
diagram 
$$\xymatrix{
C\drlowertwocell<0>{<-2>}\rto & B\dto^{\text{\rm qis}}\\
& B''}$$
Apply Theorem~\ref{loc.fin} to $C\to B''$, to obtain the homotopy
commutative 
diagram 
$$\xymatrix{
C\rto\dto_{\text{\rm finite}}\drtwocell\omit& B\dto^{\text{\rm open
imm.}} \\ 
A\rto_{\text{\rm qis}} & B'}$$
Finally, use the fact that $C\to A$ is resolving, to eliminate the
homotopy.
\end{pf}

By weakening the conclusion somewhat, we can improve on the estimate
for the amplitude of the resolution $A$:

\begin{scholum}
Let $C\to B$ be a perfect morphism of quasi-finite resolving
algebras. Then for every $K$-valued augmentation $B\to K$ there exists
an open immersion $B\to B'$, such that $B\to K$ factors through $B\to
B'$, and a commutative diagram 
$$\xymatrix{
C\rto\dto_{\text{\rm finite}} & B\dto^{\text{\rm open imm.}}\\
A\rto^-{\text{\rm \'etale}} & B'}$$
where $C\to A$ is a finite resolving morphism and $A\to B'$ 
\'etale.
If $C\to B$ has perfect amplitude $N$, then $C\to A$ can be chosen to
have amplitude $\max(1,N)$.
\end{scholum}
\begin{pf}
By examining the proof of Theorem~\ref{loc.fin}, we see that the only
place we used that $n\geq2$ was to conclude that $h^0(A)\to
h^0(B_{\{g\}})$ was an isomorphism, in the proof that $A\to B_{\{g\}}$
was a quasi-isomorphism. If we drop the assumption that $n\geq2$, but
assume that $n\geq1$, we can still conclude that $L_{B_{\{g\}}/A}$ is
acyclic, which proves that $A\to B_{\{g\}}$ is \'etale.
\end{pf}

To show that the class of perfect resolving algebras is larger than 
the class of finite resolving algebras, we provide an example of a 
perfect resolving algebra which has no resolution by a finite
resolving algebra:

\begin{numex}\label{abstract}
let $X$ be a smooth affine scheme, $E$ a vector bundle over $X$ and
$s:X\to E$ a section.  To this data we associate the Koszul complex
\begin{equation}\label{k.k}
\Lambda_{A}{M}:\quad\ldots\longrightarrow
\Lambda_{A}^2{M} 
\longrightarrow {M} \longrightarrow A\,.
\end{equation}
Here $A$ is the affine coordinate ring of $X$ and $M$ is the
projective $A$-module corresponding to $E$, i.e., the 
module of sections of the dual of $E$.  The homomorphism $M\to A$
comes from the section $s$. 
We consider (\ref{k.k}) as a differential graded algebra $B$ by placing
$\Lambda^r_AM$ in degree $-r$, so that $B^{-r}=\Lambda^r_AM$.
We note that $B$ has a universal derivation $\du:B\to \Omega_B$, and
we have
\begin{equation*}
\Omega_B\otimes_Bh^0(B)=
\bigl[\ldots\longrightarrow
\Lambda_{A}^2{M} 
\longrightarrow {M} \longrightarrow \Omega_A\bigr]\otimes_A A/I\,,
\end{equation*}
where $I$ is the ideal of $A$ defined by the image of $M\to A$, in
other words the ideal defining the zero locus of $s$. The homomorphism
$M\to \Omega_A$ is defined as the composition of $M\to A$ and the
universal derivation $A\to \Omega_A$.  We note that
$\Omega_B\otimes_Bh^0(B)$ is perfect. Moreover, it admits a
determinant, defined as
\begin{align}
\det\big(\Omega_B\otimes_Bh^0(B)\big)&=\big(\det\Omega_A\otimes \det
M^{-1}\otimes\ldots\big)\otimes_AA/I\nonumber\\
&=\big(\omega_A\otimes\det(\Lambda_AM)\big)\otimes_AA/I\,.\label{om.m}
\end{align}

By Scholum~\ref{skol}, we may find a quasi-finite resolving algebra
$B'$, together with a quasi-isomorphism $B'\to B$. Then we get a
quasi-isomorphism $\Omega_{B'}\otimes_{B'}A/I\to\Omega_B\otimes_BA/I$,
and so $B'$ is a perfect resolving algebra.  If $B$ was
quasi-isomorphic to a finite resolving algebra, then (\ref{om.m}) would
have to be a free $A/I$-module of rank one.  Of course, there are many
examples where (\ref{om.m}) is non-trivial. 
\end{numex}

\begin{numex}\label{concrete}
For a concrete example of such a perfect resolving algebra,
consider the affine cubic curve $y^2=4(x^3-x)$. The projective
completion has one additional point (a flex), called $P_0$.  We
consider the projective completion as an elliptic curve $X$ with base 
point $P_0$. Let $X'=X-\{P_0\}$.

Let $P=(0,0)$ be the origin, which is a 2-division point in $X$ and
consider the line bundle $L=\O(-P)$ on $X'$. This is a non-trivial
line bundle on the affine curve $X'$.  Because $P$ is a 2-division
point, there exists a regular function $s$ on $X'$, with a double zero
at $P$, but no other zeroes in $X'$, in other words,
$\div(s)=2P-2P_0$. To be specific, we may take $s=x$.

Choose points $Q$, $R$, $Q'$ and $R'$ on $X'$ such that $Q+R=Q'+R'=P$
in the group $X$, but $\{Q,R\}\cap\{Q',R'\}=\varnothing$. Then there
exist regular functions $a$ and $b$ on $X'$, such that
$\div(a)=P+Q+R-3P_0$ and $\div(b)=P+Q'+R'-3P_0$.  To be specific, we
may choose $a=x+\frac{1}{2}y$ and $b=x-\frac{1}{2}y$, which determines
the points $Q$, $R$, $Q'$ and $R'$ as intersections of $y^2=4(x^3-x)$
with two lines of slope $\pm 2$ through the origin.

Consider the matrix 
$$M=\frac{1}{s}\begin{pmatrix}-ab & -b^2\\
a^2 & ab\end{pmatrix}\quad\in\quad
M\big(2\times2,\Gamma(X',\O)\big)\,.$$  
One checks that 
\begin{equation}\label{resoml}
\ldots\stackrel{M}{\longrightarrow}\O^2
\stackrel{M}{\longrightarrow}\O^2
\stackrel{(a\,\,b)}{\longrightarrow}L\to 0
\end{equation}
is an infinite resolution of $L$ by free modules of rank 2.

Now let us construct a perfect resolving algebra $B$. Let
$f(x,y)=y^2-4(x^3-x)$. We start with the finite resolving algebra
$A=k[x,y]\{\xi\}/(d\xi=f)$, which resolves the affine coordinate ring
$h^0(A)=k[x,y]/f$ of $X'$.

Take sequences
$\theta_1,\theta_2,\ldots$ 
and $\eta_1,\eta_2,\ldots$ of formal variables such that
$\deg\theta_i=\deg\eta_i=-i$, for all $i\in\zz_{\geq1}$. We let
$B^\natural=A[\theta,\eta]$ be the graded algebra, free over $A$ on all
$\theta_i$ and $\eta_i$.  We turn $B^\natural$ into a differential
graded algebra $B$ by defining a differential $d$ by
$d\theta_1=d\eta_1=0$ and 
\begin{align*}
d\theta_{i+1}&=\alpha\theta_i+\gamma\eta_i-\xi\theta_{i-1}\\
d\eta_{i+1}&=\beta\theta_i+\delta\eta_i-\xi\eta_{i-1}\,,
\end{align*}
for all $i\geq1$, where we use the convention $\theta_0=\eta_0=0$.  Here,
$\alpha,\beta,\gamma,\delta\in k[x,y]$ are 
polynomials such that the image of $\bigl(\begin{smallmatrix}
\alpha&\beta\\ \gamma&\delta\end{smallmatrix}\bigr)$ in $h^0(A)$ is
equal to $M$ and $\bigl(\begin{smallmatrix}
\alpha&\beta\\
\gamma&\delta\end{smallmatrix}\bigr)^2=f(x,y)\bigl(\begin{smallmatrix} 
1&0\\ 0&1\end{smallmatrix}\bigr)$.  In fact, we may take 
$$\begin{pmatrix} \alpha&\beta\\ \gamma&\delta \end{pmatrix}=
\begin{pmatrix} x^2-x-1&y-(x^2+x-1)\\ y+(x^2+x-1)&-(x^2-x-1)
\end{pmatrix}$$
Then $B$ is a perfect resolving algebra:  the complex
$\Omega_B\otimes_Bh^0(B)=\Omega_B\otimes _Bh^0(A)$ is the direct sum
of the resolution 
(\ref{resoml}) of $L$ and the two term resolution
$$h^0(A)\Dd \xi\longrightarrow h^0(A)\Dd x\oplus h^0(A)\Dd y$$
of $\Omega_{h^0(A)}$. Thus $\Omega_B\otimes_Bh^0(B)$ is perfect, of
amplitude contained in $[-1,0]$ and $B$ is perfect of amplitude 
1. On the other hand, $B$ cannot be quasi-isomorphic to a finite
resolving algebra, because 
$$\det\big(\Omega_B\otimes_B h^0(A)\big)=L^{-1}\,$$
is not trivial.  

This is the special case of Example~\ref{abstract}, where $E$ is a
line bundle and $s$ is the zero section.
\end{numex} 

\subsection{Compatibility with limits over truncations}\label{cwlot}

The second fundamental result on perfect resolving algebras expresses
a certain compatibility between the derivations of a perfect resolving
algebra $B$ and its truncations $B_{(n)}$. Here $B_{(n)}$, for
$n\geq-1$, denotes the differential graded subalgebra generated by
$B^{\geq-n}$. Every truncation $B_{(n)}$ is a finite resolving algebra
of amplitude $n$. 

More generally, if $C\to B$ is a quasi-finite resolving morphism, then
we define $B_{(n)}$ to be the differential graded subalgebra of $B$
generated by all of $C$ and $B^{\geq-n}$.  Then every $C\to B_{(n)}$
if a finite resolving morphism of amplitude $n$.  We hope that it is
always clear from the context, which definition of $B_{(n)}$ applies.

\begin{them}\label{plim}
Let $C\to B$ be a perfect resolving morphism of differential graded
algebras. Then for every $r$ we have
$${\projectlim_n}^1 h^r\Deru_C(B_{(n)},A)=0$$
and 
$$\projectlim_n
h^r\Deru_C(B_{(n)},A)=h^r\Deru_C(B,A)$$
for all quasi-finite resolving algebras $A$ endowed with a morphism
$B\to A$. 
\end{them}
\begin{pf}
Assume that $C\to B$ has perfect amplitude $N$. 
Consider the quasi-finite resolving morphism $B_{(n)}\to B$, for
$n\geq N$. By Remark~\ref{formnp}, we know that $B_{(n)}\to B$ has
perfect amplitude $n+1$. By Proposition~\ref{ss.th}, the spectral
sequence~(\ref{ss.th.t}), computing $h^r\Deru_{B_{(n)}}(B,A)$, has
exactly one non-zero column, the column $p=n+1$. 

Now consider the canonical homomorphism
$\Omega_{B/B_{(n)}}\to\Omega_{B/B_{(n+1)}}$ of differential graded
$B$-modules.  It induces a homomorphism of spectral
sequences~(\ref{ss.th.t}). Since both spectral sequences have only one
non-zero column, but this column is off by 1, the homomorphism of
spectral sequences vanishes. We conclude that the canonical
homomorphism
$$h^r\Deru_{B_{(n+1)}}(B,A)\longrightarrow h^r\Deru_{B_{(n)}}(B,A)$$
is zero.

Next we note that we have a morphisms of short exact sequences of
differential graded $A$-modules
$$\xymatrix{
0\rto & {\Deru_{B_{(n+1)}}(B,A)} \rto\dto & {\Deru_C(B,A)} \rto\dto &
{\Deru_C(B_{(n+1)},A)}\rto\dto & 0\\
0\rto & {\Deru_{B_{(n)}}(B,A)} \rto & {\Deru_C(B,A)} \rto &
{\Deru_C(B_{(n)},A)}\rto & 0}$$
We consider the following extract from the induced morphism of long
exact cohomology sequences:
$$\xymatrix{
{}\rto & h^r\big(\Deru_{B_{(n+1)}}(B,A)\big)\dto\rto &
h^r\big(\Deru_C(B,A)\big)\dto\rto & \\
{}\rto & h^r\big(\Deru_{B_{(n)}}(B,A)\big)\rto &
h^r\big(\Deru_C(B,A)\big)\rto & }$$
Since the left vertical arrow is zero, and the right vertical arrow is
bijective, we conclude that the upper horizontal map is zero.  Thus
the whole long exact sequence on the $n+1$ level breaks up into short
exact sequences
$$\ses{h^r\Deru_C(B,A)}{}{h^r\Deru_C(B_{(n+1)},A)}{
}{h^{r+1}\Deru_{B_{(n+1)}}(B,A)}\,.$$ 
Letting $n+1$ vary, we get an inverse system of short exact
sequences. Passing to the limit, we get the following six term exact
sequence
\begin{multline*}
0\longrightarrow\projectlim h^r\Deru_C(B,A)\longrightarrow \projectlim
h^r\Deru_C(B_{(n)},A) \longrightarrow \\
\longrightarrow
\projectlim h^{r+1}
\Deru_{B_{(n)}}(B,A)\longrightarrow
{\projectlim}^1 h^r\Deru_C(B,A)\longrightarrow\\
\longrightarrow {\projectlim}^1
h^r\Deru_C(B_{(n)},A) \longrightarrow {\projectlim}^1 h^{r+1}
\Deru_{B_{(n)}}(B,A)\longrightarrow0
\end{multline*}
Since the connecting morphisms of the third projective system are zero
for sufficiently large $n$, the associated $\projectlim$ and
${\projectlim}^1$ vanish.  All connecting morphisms in the first
projective system are bijective, so the associated ${\projectlim}^1$
vanishes.  In other words, the third, fourth and sixth term in our six
term sequence vanish. This implies the theorem.
\end{pf}

\begin{scholum}\label{twostep}
The inverse system of groups $h^r\Deru_C(B_{(n)},A)$ has the property
that if an element lifts one step, then it lifts two steps. 
\end{scholum}

\begin{rmk}
If an inverse system of groups has this property, then its
${\projectlim}^1$ vanishes.
\end{rmk}

\Section{Linearization of homotopy groups}\label{sec.linear}

Occasionally, it seems unnatural to work with negative indices for
the cohomology spaces of certain differential graded modules.  We thus
adopt the usual notation $h_i=h^{-i}$ for lowering indices.

\subsection{Preliminaries}

We start with three fundamental lemmas.

\begin{lem}\label{der.diff}
Let $C\to B\to A$ be morphisms of differential graded algebras. 
Let $D\in\Deru_C(B,A)$ and $\epsilon\in A$ be homogeneous elements of
complementary degrees, i.e., such that $\deg D+\deg\epsilon=0$. Assume
$\epsilon^2=0$. Then
$h=\phi+\epsilon D$ is a morphism of graded algebras, where $\phi:B\to
A$ is the structure morphism.  
If $D$ and $\epsilon$ are cocycles, then $h$ is a
morphism of differential graded algebras.\qed
\end{lem}
%\begin{pf}
%This is a fun little excercise, or really old, depending on your point
%of view.
%\end{pf}

\begin{lem}\label{basic.sol}
Let $A$ be a differential graded algebra and $\eta\in
A\otimes\Omega_\ell$ a 
homogeneous element, where $\ell\geq1$. Let
$\Lambda\subset\Delta^\ell$ be a horn. If 
$\eta\resto_\Lambda=0$ and $d\eta=0$, then there exists  
$\theta\in A\otimes \Omega_\ell$ such that
$\theta\resto_\Lambda=0$ and 
$\eta=d\theta$.
\end{lem}
\begin{pf}
Without loss of generality, let $\Lambda\hookrightarrow\Delta^\ell$ be
the horn defined by $\Lambda=\{t_1=0\}\cup\ldots\cup\{t_\ell=0\}$.

Note that $A\otimes\Omega_\ell\to A;\omega\mapsto\omega(0)$ is a
quasi-isomorphism, by the algebraic de Rham theorem.
Thus we can certainly find ${\theta}_0\in A\otimes\Omega_\ell$ such that
$d{\theta}_0=\eta$. We
define $\theta_i$ inductively by
$\theta_i=\theta_{i-1}-\theta_{i-1}\resto_{t_i=0}$. Then
$\theta=\theta_\ell$ fits our requirements.
\end{pf}

\begin{lem}\label{extrapolate}
Let $A$ be a differential graded algebra and let  $\psi\in
A\otimes\Omega_{\ell-1}$ be a homogeneous element, where
$\ell\geq1$. Assume that 
$\psi\resto_{\del\Delta^{\ell-1}}=0$. Then there exists $\Psi\in
A\otimes\Omega_{\ell}$ such that $\del_{\ell}\Psi=\psi$ and
$\Psi\resto\Lambda=0$, where $\Lambda\subset\Delta^\ell$ is the horn
opposite to the face
$\del^{\ell}(\Delta^{\ell-1})\subset\Delta^{\ell}$. 
\end{lem}
\begin{pf}
Note that the face map $\del^\ell:\Delta^{\ell-1}\to\Delta^\ell$ is
given in inhomogeneous coordinates by $(t_1,\ldots,t_{\ell-1})\mapsto
(t_1,\ldots, t_{\ell-1},0)$. 
Define 
$$\Psi(t_1,\ldots,t_{\ell})=
(1-t_{\ell})^{N+1}
\psi(\textstyle \frac{t_1}{1-t_{\ell}},\ldots,
\frac{t_{\ell-1}}{1-t_\ell})\,,$$   
where $N$ is large enough such that $(1-t_{\ell})^{N}
\psi(\textstyle \frac{t_1}{1-t_{\ell}},\ldots,
\frac{t_{\ell-1}}{1-t_\ell})$ has no denominators.
\end{pf}

\begin{rmk}
The latter two lemmas can be interpreted in terms of the differential
graded algebra $B_r=k[x,\xi]/(d\xi=x)$, which is quasi-free on the basis
$(x,\xi)$, where $\deg x=r$ and $\deg\xi=r-1$, and where the
differential is defined by $d\xi=x$, $dx=0$.  Lemma~\ref{basic.sol}
expresses the fact that $\shom(B_r,A)\to\shom(k[x],A)$ is a
fibration, which also follows from Corollary~\ref{qfreefib}.
Lemma~\ref{extrapolate} says that 
$\pi_{\ell-1}\shom(B_r,A)=0$. This also follows from the fact that
$\shom(B_r,A)\to \shom(k,A)=\ast$ is a trivial fibration, since $k\to
B_r$ is a quasi-isomorphism, which follows most easily from
Corollary~\ref{qiscondition}. 
\end{rmk}

Let $A$ be a differential graded algebra.

Consider the simplicial differential graded algebra
$A\otimes\Omega\lcom$. Let $N(A\otimes\Omega\lcom)$ be the associated
normalized chain complex.  Thus $N(A\otimes\Omega_\ell)\subset
A\otimes\Omega_\ell$ is the subcomplex defined by
$$\omega\in N(A\otimes\Omega_\ell)\quad\Longleftrightarrow
\quad\text{for all $i=0,\ldots,\ell-1$ we have $\del_i\omega=0$,}$$ 
where $\del_i\omega=(\del^i)\upst\omega$ is restriction via the $i$-th
inclusion map 
$\del^i:\Delta^{\ell-1}\longrightarrow\Delta^\ell$.
The boundary map of the normalized chain complex is defined to be
$$\tilde{\del}_\ell=(-1)^\ell\del_\ell:
N(A\otimes\Omega_\ell)\longrightarrow  N(A\otimes\Omega_{\ell-1})\,.$$ 
Written in such a way, $\tilde{\del}_\ell$ is a morphism of complexes,
i.e., it commutes with the coboundary map $d$ of
$N(A\otimes\Omega\lcom)$.  Since we would like to think of
$\tilde{\del}$ as being of degree $-1$,  we have to change the sign of
$d$ on $N(A\otimes\Omega_\ell)$ for odd $\ell$.  We call the result
$\tilde{d}$. Thus we have
$$\tilde{d}\omega=(-1)^\ell d\omega,\quad\text{for $\omega\in
N(A\otimes\Omega_\ell)$.}$$

Now we truncate each of the cochain complexes
$N(A\otimes\Omega_\ell)$ at $r$. We obtain a chain complex of cochain
complexes, whose boundary maps are given by
$$\tilde{\del}_\ell=(-1)^\ell\del_\ell:\big(\tau_{\leq
r}N(A\otimes\Omega_\ell),\tilde{d}\,\big)\longrightarrow \big(
\tau_{\leq r}N(A\otimes\Omega_{\ell-1}),\tilde{d}\,\big)\,.$$
Thus we have defined a double complex, which we write in the third
quadrant, with $\tilde{d}$ vertical and $\tilde{\del}$ horizontal:
$$\xymatrix{
\ldots\ar[r]^-{-\del_3} & Z^rN(A\otimes\Omega_2)\ar[r]^-{\del_2} & Z^r
N(A\otimes\Omega_1)\ar[r]^-{-\del_1} & Z^r A\\
\ldots\ar[r]^-{-\del_3} & N(A\otimes\Omega_2)^{r-1}\uto_d\ar[r]^-{\del_2} & 
N(A\otimes\Omega_1)^{r-1}\uto_{-d}\ar[r]^-{-\del_1} & A^{r-1}\uto_d\\
\ldots\ar[r]^-{-\del_3} & N(A\otimes\Omega_2)^{r-2}\uto_d\ar[r]^-{\del_2} & 
N(A\otimes\Omega_1)^{r-2}\uto_{-d}\ar[r]^-{-\del_1} & A^{r-2}\uto_d\\
&{\vdots}\uto_d & {\vdots}\uto_{-d} & {\vdots}\uto_d}$$
By Lemma~\ref{basic.sol}, vertical cohomology vanishes everywhere,
except in the last column. As for the horizontal cohomology, it
vanishes everywhere except for in the first row, by
Lemma~\ref{extrapolate}.
Therefore, our double complex induces an isomorphism of $k$-vector
spaces
\begin{equation}\label{th.es}
h^{r-\ell}(A)\longiso
h_{\ell}\big(Z^rN(A\otimes\Omega\lcom),\tilde{\del}\,\big) \,.
\end{equation}

We will fix notation: 
\begin{equation}\label{om}
\omega_\ell=\ell!\,dt_1\wedge\ldots\wedge dt_\ell\,.
\end{equation}
Note the properties $d\omega_\ell=0$ and
$\omega_\ell\resto_{\del\Delta^\ell}=0$. Also, $\omega_0=1$.

\begin{prop}\label{prep.a}
Every element of
$h_{\ell}\big(Z^rN(A\otimes\Omega\lcom),\tilde{\del}\,\big)$ may be
represented as $\omega_\ell a$, where $a\in Z^{r-\ell}(A)$. Changing
$a$ by a coboundary only changes $\omega_\ell a$ by a boundary.

The canonical isomorphism (\ref{th.es}) is given by 
$$a\longmapsto (-1)^{\frac{1}{2}\ell(\ell-1)}\omega_\ell a\,.$$
\end{prop}
\begin{pf}
We shall introduce differential forms (for $\ell\geq1$)
$$\tau_\ell=-(\ell-1)!\sum_{i=1}^\ell(-1)^i
t_i\,dt_1\wedge\ldots\wedge\widehat{dt_i}\wedge\ldots\wedge
dt_\ell\,$$
Note the following properties of $\tau_\ell$:

(i) $d\tau_\ell=\omega_\ell$,

(ii) for all $i=1,\ldots,\ell$ we have $\tau_\ell\resto_{t_i=0}=0$,
i.e., $\del_i\tau_\ell=0$,

(iii) $\del_0\tau_\ell=\omega_{\ell-1}$.

\noindent%
To see (iii), recall that
$\del^0(t_1,\ldots,t_{\ell-1})=(1-\sum_{i=1}^{\ell-1}
t_i,t_1,\ldots,t_{\ell-1})$.

For the present proof, it is convenient to consider the following
variation:
$$\sigma_\ell=\tau_\ell+(-1)^\ell\omega_{\ell-1}\,.$$
These forms have the three properties:

(i) $d\sigma_\ell=\omega_\ell$,

(ii) for all $i=0,\ldots,\ell-1$ we have $\del_i\sigma_\ell=0$,

(iii) $\del_\ell\sigma_\ell=(-1)^\ell\omega_{\ell-1}$, so that
$\tilde{\del}_\ell\sigma_\ell= \omega_{\ell-1}$.

Let $a\in Z^{r-\ell}(A)$ and consider the sum
\begin{equation}\label{zz}
\sum_{i=1}^\ell(-1)^{\frac{1}{2}i(i+1)}\sigma_ia\,.
\end{equation}
Note that $\sigma_ia\in N(A\otimes\Omega_i)^{r-\ell+i-1}$, so that
(\ref{zz}) is an element of total (cochain) degree
$r-\ell-1$ in our double complex. Applying the total coboundary
$\tilde{d}+\tilde{\del}$ to 
(\ref{zz}) we get
\begin{align*}
(\tilde{d}+\tilde{\del})
\sum_{i=1}^\ell(-1)^{\frac{1}{2}i(i+1)}\sigma_ia &=
\sum_{i=1}^\ell(-1)^{\frac{1}{2}i(i+1)+i}d\sigma_ia +
\sum_{i=1}^\ell(-1)^{\frac{1}{2}i(i+1)}\tilde{\del}_i\sigma_ia\\
&=\sum_{i=1}^\ell(-1)^{\frac{1}{2}i(i+1)+i}\omega_ia +
\sum_{i=1}^\ell(-1)^{\frac{1}{2}i(i+1)}\omega_{i-1}a\\
&=\sum_{i=1}^\ell(-1)^{\frac{1}{2}i(i-1)}\omega_ia +
\sum_{i=0}^{\ell-1}(-1)^{\frac{1}{2}(i+1)(i+2)}\omega_{i}a\\
&=(-1)^{\frac{1}{2}\ell(\ell-1)}\omega_\ell a-a\,.
\end{align*}
This proves the formula of the proposition.  It also proves  the first
claim, because (\ref{th.es}) is surjective. For the second claim let
$b\in A^{r-\ell-1}$. Then 
$$\psi=\sigma_{\ell+1}db+(-1)^\ell\omega_{\ell+1}b$$
is an element of $Z^rN(A\otimes\Omega_{\ell+1})$
and $\tilde{\del}\psi=\omega_\ell db$.
\end{pf}

\subsection{Linearization of homotopy groups}\label{sec.lin}

Let $C\to B$ be a fixed resolving morphism of resolving algebras over
$k$.
Let
$A$ be an arbitrary differential graded algebra and denote by  $P:B\to A$ a
fixed morphism. Denote by $P$ also the restriction to $C$.  We use $P$
as a base point for the spaces $\shom(B,A)$ and $\shom(C,A)$ and for
the fiber $\shom_C(B,A)$ of the fibration $\shom(B,A)\to\shom(C,A)$. 

Consider, for every $\ell\geq1$, the map
\begin{align}\label{def.xil}
\Xi_\ell:h^{-\ell}\Deru_C(B,A)&\longrightarrow
\pi_\ell\shom_C(B,A)\nonumber\\*
D&\longmapsto P+(-1)^{\frac{1}{2}\ell(\ell-1)}\omega_\ell D\,.
\end{align}
Recall that $\omega_\ell$ has been defined in~(\ref{om}).
By Lemma~\ref{der.diff}, for given $D\in Z^{-\ell}\Deru_C(B,A)$, the
map 
$$h=P+(-1)^{\frac{1}{2}\ell(\ell-1)}\omega_\ell D$$
is a morphism of
differential graded $C$-algebras $B\to A\otimes\Omega_\ell$, hence an
$\ell$-simplex in $\shom_C(B,A)$. Since
$\omega_\ell\resto_{\del\Delta^\ell}=0$, it defines, in fact, an
element of 
$\pi_\ell\shom_C(B,A)$. 

Let us check that the homotopy class of $h$ depends only on the
cohomology class of $D$.  So let $D'$ be another element of
$Z^{-\ell}\Deru_C(B,A)$, differing from $D$ by a coboundary and let
$h'=P+(-1)^{\frac{1}{2}\ell(\ell-1)}\omega_\ell D'$. Thus there
exists $E\in\Deru^{-\ell-1}_C(B,A)$
such that 
$$D'=D+dE\,.$$

Then one checks that
$$\Phi=P+(-1)^{\frac{1}{2}\ell(\ell-1)}\omega_\ell\big((1-s)D+sD'\big)+
(-1)^{\frac{1}{2}\ell(\ell-1)}\omega_\ell 
ds \,E$$
defines a homotopy from $h$ to $h'$, i.e., an element 
$$\Phi\in\Hom\big(I\times\Delta^\ell,\shom_C(B,A)\big)\,,$$
such that
$\Phi\resto_{s=0}=h$, $\Phi\resto_{s=1}=h'$ and
$\Phi\resto_{I\times\del \Delta^\ell}=P$.  Here $s$ is the coordinate
on the `interval' $I=\aaa^1$.

Thus the map $\Xi_\ell$ is well-defined, for all $\ell\geq1$.

\begin{numrmk}[Naturality]\label{nature.xi}
(i) Let $C'\to B'$ be a another resolving morphism of
resolving algebras. Assume, moreover, given a
commutative diagram of differential graded algebras
$$\xymatrix@=15pt{
C'\dto\rto & B'\dto\\
C\rto & B\,.}$$
Then we have, for all $\ell\geq1$,  an induced commutative diagram
$$\xymatrix{
h^{-\ell}\Deru_C(B,A) \dto_{\Xi_\ell}\rto &
h^{-\ell}\Deru_{C'}(B',A) \dto^{\Xi_\ell}\\
\pi_{\ell}\shom_C(B,A)\rto&\pi_{\ell}\shom_{C'}(B',A)\,.}$$

(ii) Let $A\to A'$ be an arbitrary morphism of differential graded
algebras. Then we have an induced commutative diagram
$$\xymatrix{
h^{-\ell}\Deru_C(B,A) \dto_{\Xi_\ell}\rto &
h^{-\ell}\Deru_{C}(B,A') \dto^{\Xi_\ell}\\
\pi_{\ell}\shom_C(B,A)\rto&\pi_{\ell}\shom_{C}(B,A')\,.}$$

(iii) As a special case of~(ii), applied to
$\del_0,\del_1:A\otimes\Omega_1\toto A$, we get the following
compatibility of $\Xi_\ell$ with change of base point: let $A$ be a
differential graded $C$-algebra (whose cohomology is concentrated in
non-positive degrees). Let $P,Q:B\to A$ be two base points for
$\shom_C(B,A)$, giving two different $B$-modules structures on $A$,
denoted ${}_PA$, ${}_QA$. Then every path $h:P\to Q$ in
$\shom_C(B,A)$ gives rise to a commutative diagram
$$\xymatrix{
{h^{-\ell}\Deru_C(B,{}_PA)} \dto_{\Xi_\ell}\rto^{h\lst} & 
{h^{-\ell}\Deru_C(B,{}_QA)} \dto^{\Xi_\ell}\\
{\pi_\ell\big(\shom_C(B,A),P\big)}\rto^{{}^h(\argument)} &
{\pi_\ell\big(\shom_C(B,A),Q\big)}}$$
where $h\lst$ denotes the change of base point map of
Definition~\ref{can-iso}.
\end{numrmk}

\begin{lem}[Homomorphism]\label{Homo}
Let $C\to B\to A$ be as above.

(i) For $\ell\geq2$, the map $\Xi_\ell$ is a group homomorphism.  

(ii) If there
exists a basis $(x_\nu)$ for $B$ over $C$,
such that $dx_\nu\in C$, for all $\nu$, then this is also true for
$\ell=1$.

(iii) More generally, suppose that $B'\subset B$ is a subalgebra
containing $C$, such that $C\to B'$ and $B'\to B$ are resolving
and there exists a basis $(x_\nu)$ for $B$
over
$B'$, such that $dx_\nu\in B'$, for all $\nu$. Let $D\in
h^{-1}\Deru_C(B,A)$ and $D'\in h^{-1}\Deru_{B'}(B,A)$. Then we have
$\Xi_\ell(D+D')=\Xi_\ell(D)\ast\Xi_\ell(D')$. 
\end{lem}
\begin{pf}
First assume that $\ell\geq2$.
Let $D$, $D'$ be two elements of $Z^{-\ell}\Deru_C(B,A)$.
Let $h$, $h'$
and $g$ be the images of $D$, $D'$ and $D+D'$ under $\Xi_\ell$. An
$(\ell+1)$-simplex in $\shom_C(B,A)$ showing that $h\ast h'=g$ is given
by
\begin{align*}
\Phi=P+&\epsilon(dt_1\wedge\ldots\wedge dt_{\ell-2}\wedge
dt_\ell\wedge dt_{\ell+1}+dt_1\wedge\ldots\wedge dt_{\ell-1}\wedge 
dt_{\ell+1})\,D \\ 
&\epsilon(dt_1\wedge\ldots\wedge dt_{\ell-1}\wedge 
dt_{\ell+1}+dt_1\wedge\ldots\wedge dt_\ell)\,D', 
\end{align*}   
where $\epsilon=(-1)^{\frac{1}{2}\ell(\ell-1)}\ell!$ is the necessary
multiplier. Note that we need $\ell\geq2$ to apply Lemma~\ref{der.diff}.
We see
directly from the definition that $\Phi\resto\{t_{\ell-1}=0\}=h$, that
$\Phi\resto\{t_{\ell+1}=0\}=h'$, that $\Phi\resto\{t_{\ell}=0\}=g$ and
that $\Phi$ restricted to all the other faces of $\Delta^{\ell+1}$ is
equal to $P$.

Now let us consider the case $\ell=1$. Since (iii) implies (ii), let
us prove (iii). So 
assume given $C\to B'\to B$. We define $\Phi$
to be the unique morphism of graded algebras $\Phi:B\to
A\otimes\Omega_2$ such that $\Phi\resto B'=P+dt_2\,D$ and 
$$\Phi(x_\nu)=P(x_\nu)+ dt_2\,D(x_\nu)+(dt_2+dt_1)D'(x_\nu)\,,$$
for all $\nu$.  One checks that $\Phi$ respects the differential, so
that $\Phi$ is a 2-simplex in $\shom_C(B,A)$. Moreover,
$\Phi\resto_{t_1=1-t\atop t_2=t\phantom{-1}}=h$,
$\Phi\resto_{t_1=0}=g$ and $\Phi\resto_{t_2=0}=h'$, so that, indeed,
$h\ast h'=g$.
\end{pf}

\begin{numrmk}[Degree zero]
Suppose $B$ admits a basis  $(x_i)$ over $C$ 
such that $dx_i\in C$, for all $i$. Then we may also define 
$$\Xi_0:h^0\Deru_C(B,A)\longrightarrow\pi_0\shom_C(B,A)\,.$$
Given a $C$-derivation $D:B\to A$, we map $D$ to the unique morphism
of graded algebras $h:B\to A$ such that $h\resto C=P$ and
$h(x_i)=P(x_i) + D(x_i)$. To check that $h$ respects the differential,
it suffices to prove that $h(dx_i)=d\big(h(x_i)\big)$, for all $i$,
which is easily done using that $h(dx_i)=P(dx_i)$ and $D(dx_i)=0$,
which follows from our assumption that $dx_i\in C$.  If $D'=D+dE$ is
in the same cohomology class as $D$, then the image of $D$, which we
call $h'$, is homotopic to $h$ via the homotopy $\Phi$ defined by
$\Phi\resto C=P$ and $\Phi(x_i)=P(x_i)+(1-s)D(x_i)+s D'(x_i)+ds\,
E(x_i)$. Thus $\Xi_0$ is indeed well-defined. 

The question of whether or not $\Xi_0$ depends on the choice of
generators $(x_i)$ is slightly more subtle.  The most common reason
why a set of generators should satisfy $dx_i\in C$ is that they all
have the same degree $r=\deg x_i$, for all $i$.  So let us suppose
that this is the case, and that the total number of generators $x_i$
if finite.  

In the case $r<0$ it is easy to see that $\Xi_0$ is
independent of the choice of the generators: if $(y_j)$ is another
family of generators, then $x_i=\sum_jc_{ij}y_j$, for a family of
elements $c_{ij}$ of degree zero, and hence necessarily contained in
$C$. So if $h'$ is defined by $h'(y_j)=P(y_j) +D(y_j)$, then
$h'(x_i)=h'(\sum_j
c_{ij}y_j)=\sum_jP(c_{ij})\big(P(y_j)+D(y_j)\big)=P(\sum_jc_{ij}y_j)+
D(\sum_jc_{ij} y_j)=h(x_i)$ and so $h'=h$.

On the other hand, in the case $r=0$, the map $\Xi_0$ depends on the
choice of generators. For example, let $C=k$ and
$B=k[x_1,x_2]$. Then
another set of free generators for $B$ is given by $y_1=x_1$ and
$y_2=x_2+x_1^2$. The fact that the change of coordinates is not linear
is responsible for the fact that $h$ and $h'$ will be different. This
is particularly easy to see if $A$ is also concentrated in degree
zero, because then $\pi_0\shom(B,A)=\Hom(B,A)$, the set of $k$-algebra
morphisms.
\end{numrmk}

The reason for the sign in the definition of $\Xi_\ell$ becomes clear
from the following fundamental lemma:

\begin{lem}[Main]\label{Main}
Suppose that $B$ has a basis over $C$ consisting of one homogeneous
generator $x$. Let $r=\deg x$. Then for all 
$\ell\geq0$ the diagram
$$\xymatrix{
h^{-\ell}\Deru_C(B,A)\dto_{\Xi_\ell}\ar[rrr]^{D\mapsto D(x)} &&&
h^{r-\ell}(A)\dto^{\text{\upshape can.\ hom.\ (\ref{th.es})}}\\
\pi_\ell\shom_C(B,A)\ar[rrr]_{h\mapsto h(x)-P(x)} &&& h_\ell\big(
Z^rN(A\otimes\Omega\lcom)\big)}$$
commutes. Moreover, all the maps are bijections of pointed sets, hence
group isomorphisms for $\ell\geq1$.
\end{lem}
\begin{pf}
It is straightforward to check that the upper horizontal map is
well-defined and an isomorphism of $h^0(A)$-modules. As for the lower
horizontal map, note that we have in fact
an isomorphism of pointed simplicial sets
\begin{align*}
\shom_C(B,A)&\longrightarrow Z^r(A\otimes\Omega\lcom)\\*
h&\longmapsto h(x)-P(x)\,.
\end{align*}
Moreover, the pointed simplicial set $Z^r(A\otimes\Omega\lcom)$ 
is a simplicial $k$-vector space and so its homotopy groups are equal
to the homology groups of the associated normalized chain complex. So
we have
$$\pi_\ell\big(Z^r(A\otimes\Omega\lcom)\big)= h_\ell\big(N
Z^r(A\otimes\Omega\lcom)
\big)=h_\ell\big(Z^rN(A\otimes\Omega\lcom)\big)\,.$$
So, indeed, the lower horizontal map is an isomorphism of groups
(or pointed sets, if $\ell=0$).  Commutativity of the diagram follows
directly from the definitions.
\end{pf}

\begin{cor}\label{soepf}
Suppose that $B$ has a $C$-basis  $(x_\nu)$, where all $x_\nu$ have
the same degree $r$. 
Then $\Xi_\ell$ induces canonical isomorphisms
$$h^{-\ell}\Deru_C(B,A)=\pi_\ell\shom_C(B,A)=h^{r-\ell}(A)^{\#\nu}\,,$$
for all $\ell\geq0$.
\end{cor}

Let $C\to B'\to B$ be resolving morphisms of resolving algebras.
Consider the fibration 
\begin{equation}\label{fibr.der}
\shom_C(B,A)\longrightarrow\shom_C(B',A)\,,
\end{equation}
whose fiber over $P$ is $\shom_{B'}(B,A)$. 
There is also a short exact sequence of complexes of $k$-vector spaces
\begin{equation}\label{ses.der}
\ses{\Deru_{B'}(B,A)}{}{\Deru_C(B,A)}{}{\Deru_C(B',A)}\,.
\end{equation}

\begin{prop}[Boundary]\label{boundary}
Assume, moreover, that there exists a basis
$(x_\nu)$ for $B$ 
over $B'$ such that $dx_\nu\in B'$, for all $\nu$. Then for all
$\ell\geq1$, the diagram 
$$\xymatrix{
h^{-\ell}\Deru_C(B',A)\dto_{\Xi_\ell}\rto^\delta&
h^{1-\ell}\Deru_{B'}(B,A)\dto^{\Xi_{\ell-1}}\\
\pi_\ell\shom_C(B',A)\rto^\del &\pi_{\ell-1}\shom_{B'}(B,A)}$$
commutes.
Here $\delta$ is the boundary map of the long exact 
cohomology sequence associated to~(\ref{ses.der}) and $\del$ is the
boundary map of the long exact homotopy sequence associated
to~(\ref{fibr.der}). 
\end{prop}
\begin{pf}
Consider $D\in Z^{-\ell}\Deru_C(B',A)$.  Lift $D$ to an internal
derivation $\tilde{D}:B\to A$ by setting $\tilde{D}(x_\nu)=0$, for all
$\nu$. Then $\delta D$ is equal to 
$d\tilde{D}$. Thus $\delta D$ is given by 
$$(\delta D)(x_\nu)=(d\tilde{D})(x_\nu)=(-1)^{\ell-1} D(dx_\nu)\,.$$
Applying $\Xi_{\ell-1}$ we obtain
the $\ell-1$ simplex 
$$g=P+(-1)^{\frac{1}{2}(\ell-1)(\ell-2)}\omega_{\ell-1}\delta{D}\,\in 
\pi_{\ell-1}\shom_{B'}(B,A)\,.$$
Note that 
\begin{equation}\label{gxnu}
g(x_\nu)=P(x_\nu)+
(-1)^{\frac{1}{2}\ell(\ell-1)}\omega_{\ell-1}D(dx_\nu)\,. 
\end{equation}

If, on the other hand, we apply $\Xi_\ell$ to $D$, we obtain the
$\ell$-simplex 
$$h=P+(-1)^{\frac{1}{2}\ell(\ell-1)}\omega_\ell D\,
\in\pi_\ell\shom_C(B',A)\,.$$  
We lift $h$ to $h':B\to A\otimes\Omega_\ell$ by 
$$h'(x_\nu)=P(x_\nu)+(-1)^{\frac{1}{2}\ell(\ell-1)}\tau_\ell
D(dx_\nu)\,,$$
where $\tau_\ell$ is the $\ell-1$ form defined in the proof of
Proposition~\ref{prep.a}. 
Note that $d\big(h'(x_\nu)\big)=h(dx_\nu)$, so that this formula defines a
morphism of differential graded algebras $B\to
A\otimes\Omega_\ell$. Note also that $h'\resto\{t_i=0\}=P$, for all
$i=1,\ldots,\ell$. Thus the  image of $h$ under $\del$ is equal to
$\del_0(h')$. We obtain $\del h\resto B'=P$ and 
\begin{multline*}
\del h(x_\nu)=\del_0 h'(x_\nu) \\ =P(x_\nu)+
(-1)^{\frac{1}{2}\ell(\ell-1)}\del_0\tau_\ell D(dx_\nu)
=P(x_\nu)+
(-1)^{\frac{1}{2}\ell(\ell-1)}\omega_{\ell-1} D(dx_\nu)\,,
\end{multline*}
which agrees with (\ref{gxnu}), so we see that $\del h=g$.
\end{pf}

\begin{them}\label{Bij}
Assume that $B$ is finite as a resolving algebra over $C$. Then for
$\ell\geq1$, the map 
$$\Xi_\ell:h^{-\ell}\Deru_C(B,A)\longrightarrow\pi_\ell\shom_C(B,A)$$
is bijective (so for $\ell\geq2$ an isomorphism). 
\end{them}
\begin{pf}
Induction on the number $n$ of elements in a basis of $B$ over $C$. If
this number 
is zero, then $B=C$ and so the claim is trivial. Otherwise, there
exists a differential graded subalgebra $B'\subset B$, such that $B$
has a basis over $B'$ consisting of  one homogeneous generator $x$,
and $B'$ has a $C$-basis with $n-1$ elements.  Note that then $dx\in B'$. 

Consider the fibration (\ref{fibr.der}) 
and the short exact sequence of complexes of $k$-vector
spaces~(\ref{ses.der}). 
We have the associated long exact sequences of homotopy groups and
cohomology spaces. The maps $\Xi_\ell$ relate the two:
\begin{multline*}
\xymatrix{
\ldots\rto &
h^{-\ell}\Deru_{B'}(B,A)\dto_{\Xi_\ell}\rto& 
h^{-\ell}\Deru_C(B,A)\dto^{\Xi_\ell}\rto&\\
\ldots\rto&
\pi_\ell\shom_{B'}(B,A)\rto&
\pi_\ell\shom_C(B,A)\rto&}\\
\xymatrix{
{}\rto&
h^{-\ell}\Deru_C(B',A)\dto_{\Xi_\ell}\rto^\delta&
h^{1-\ell}\Deru_{B'}(B,A)\dto^{\Xi_{\ell-1}}\rto&
\ldots\\ 
{}\rto&
\pi_\ell\shom_C(B',A)\rto^\del&
\pi_{\ell-1}\shom_{B'}(B,A)\rto&\ldots}
\end{multline*}
By naturality and Proposition~\ref{boundary} all squares which are
defined commute.
Let us first consider the part up to
$\Xi_1:h^{-1}\Deru_{B'}(B,A)\to 
\pi_1\shom_{B'}(B,A)$. This is, by Lemma~\ref{Homo}, a homomorphism of
long exact sequences of abelian groups. By induction and
Lemma~\ref{Main}, all $\Xi_\ell$ adjacent to a square involving
$\delta$ and $\del$ are isomorphisms. Applying the 5-lemma proves the
theorem in the $\ell\geq2$ case.

For the $\ell=1$ case, we look at the last five terms that still have
a map $\Xi_\ell$ defined.  This part ends at $\Xi_0:h^{0}\Deru_{B'}(B,A)\to 
\pi_0\shom_{B'}(B,A)$. By induction and Lemma~\ref{Main} all five maps
are bijective, except for the one in the middle, which is
$\Xi_1:h^{-1}\Deru_C(B,A)\to\pi_{1}\shom_C(B,A)$. To show that this
map is also bijective, note that is is equivariant for the action of
$h^{-1}\Deru_{B'}(B,A)$, because of Lemma~\ref{Homo}~(iii). This is
sufficient to prove our claim, by a suitably generalized
5-lemma.
\end{pf}

Recall from Section~\ref{cwlot} the truncations $B_{(n)}$ of the
perfect resolving morphism $C\to B$.
\begin{cor}\label{greatc}
Let $C\to B$ be a perfect resolving morphism of resolving
algebras. Then for any differential graded $B$-algebra $A$, the canonical
homomorphism 
$$\pi_\ell\shom_C(B,A)\longiso \projectlim_n \pi_\ell\shom_C(B_{(n)},A)$$
is an isomorphism, for all $\ell\geq0$.  Moreover,
$$\Xi_\ell:h^{-\ell}\Deru_C(B,A)\longrightarrow\pi_\ell\shom_C(B,A)$$
is bijective, for all $\ell\geq1$.
\end{cor}
\begin{pf}
Directly from the definitions, it follows that $\shom_C(B_{(n)},A)$ is
a tower of fibrations and that
$$\shom_C(B,A)=\projectlim_n\shom_C(B_{(n)},A)\,.$$
From Theorem~\ref{Bij} and Theorem~\ref{plim} we get that
\begin{equation}\label{limone}
{\projectlim_n}^1\pi_\ell\shom_C(B_{(n)},A)=0\,,
\end{equation}
for all $\ell\geq1$. 
In fact, for the case $\ell=1$, we have to be a little careful,
because $\Xi_1$ is not a group homomorphism, and in general,
${\projectlim}^1$ depends on the group structure, not only on the
inverse system structure.  But in Scholum~\ref{twostep} we mentioned a
property of inverse systems which $h^{-1}\Deru_C(B_{(n)},A)$ enjoys,
and which does not depend on the group structure, and is thus
inherited by $\pi_1\shom_C(B_{(n)},A)$.  As we remarked above
(following the scholum), this property is sufficient to assure that
${\projectlim}^1$ vanishes.

So the first claim follows immediately from the
Milnor exact sequence (see 
\cite{simhomthe}, Chapter~VI, Proposition~2.15).  The second claim
then follows by taking the limit over Theorem~\ref{Bij}, applied to
the various truncations $B_{(n)}$.
\end{pf}

Let us introduce the notation
$$\tilde{\omega}_\ell=(-1)^{\frac{1}{2}\ell(\ell-1)}\omega_\ell\,.$$
Let $B$ be a finite resolving algebra over $k$.  Let $A$ be an
arbitrary 
differential graded $k$-algebra and $P:B\to A$ a base point for
$\shom(B,A)$. Let $\ell\geq1$.

By the theorem,  every element of $\pi_\ell\shom(B,A)$
may be written as
\begin{equation}\label{stan}
h=P+\tilde{\omega}_\ell D\,,
\end{equation}
for a unique $D\in h^{-\ell}\Deru(B,A)$. We may call~(\ref{stan}) the
{\em standard }form of $h$. For $\ell\geq2$, we have
$$(P+\tilde{\omega}_\ell D)\ast(P+\tilde{\omega}_\ell D') = 
P+\tilde{\omega}_\ell(D+D')\,,$$
but we have no such simple formula in the case $\ell=1$.

It is natural to ask if there is nonetheless some way to describe the
composition in $\pi_1\shom(B,A)$ in terms of $h^{-1}\Deru(B,A)$.

\subsection{What can we say about $\pi_0$?}\label{wcwsa}

We will study more carefully $\pi_0\shom(B,A)$, in particular how it
behaves under change of $B$.

Let  $C\to B'\to B$ be resolving morphisms of resolving
algebras. Let $A$ be a differential
graded algebra and let us fix a morphism $P:B'\to A$, but let us not
fix any morphism $B\to 
A$.   Thus, in the fibration $\shom_C(B,A)\to\shom_C(B',A)$,
only the base is pointed, the total space is not.  We still have a
well-defined fiber $\shom_{B'}(B,A)$, but this fiber is not pointed
either.  We write the tail end of the long exact homotopy sequence as
\begin{equation}\label{tail}
\xymatrix@C=10pt{
\pi_1\shom_C(B',A)\ar@{o}[r]&\pi_0\shom_{B'}(B,A)\rto &
\pi_0\shom_C(B,A)\rto & \pi_0\shom_C(B',A)\,.}
\end{equation}
This means that $\pi_1\shom_C(B',A)$ acts on $\pi_0\shom_{B'}(B,A)$ in
such a way that two elements are in the same orbit if and only if they
map to the same element of $\pi_0\shom_C(B,A)$.  Moreover, an element
of $\pi_0\shom_C(B,A)$ lifts to $\pi_0\shom_{B'}(B,A)$ if and only if
it maps to the base point $P$ of $\pi_0\shom_C(B',A)$. 

Our first goal is an amplification of (\ref{tail}). For this we assume
that $B$ has a basis over $B'$ consisting of one element $x$ of
degree $r$, which we shall fix throughout our discussion.  Note that
$dx\in B'$.  We will also assume that $B'$ is finite over
$C$, so that we may apply Theorem~\ref{Bij}.

\begin{prop}\label{xi}
The $k$-vector space $h^r(A)$ acts transitively on the fiber of 
$$\pi_0\shom_C(B,A)\longrightarrow\pi_0\shom_C(B',A)$$
over $P$.  The stabilizer of this action is the image of the
homomorphism
\begin{align*}
\xi_P:h^{-1}\Deru_C(B',A)&\longrightarrow h^r(A)\\*
D&\longmapsto D(dx)\,.
\end{align*}
\end{prop}
\begin{pf}
Let us start by defining an action of $h^r(A)$ on
$\pi_0\shom_{B'}(B,A)$. Given $a\in h^r(A)$ and
$Q\in\pi_0\shom_{B'}(B,A)$, we define $a\ast Q\in
\pi_0\shom_{B'}(B,A)$ to the the unique element represented by the
morphism $h:B\to A$, such that $h\resto_{B'}=P$ and $h(x)=Q(x)+a$.  

To check that this is well-defined, let $a'=a+db$ and $Q':B\to A$ be
homotopic to $Q$.  Choose a homotopy $H:B\to A\otimes\Omega_1$ such
that $H\resto_{B'}=P$, $H(0)=Q$, $H(1)=Q'$. Let $h:B\to A$ satisfy
$h\resto_{B'}=P$ and $h(x)=Q(x)+a$. Let $h':B\to A$ satisfy
$h'\resto_{B'}=P$ and $h'(x)=Q'(x) +a'$.  Then define a homotopy
$H':B\to A\otimes \Omega_1$ from $h$ to $h'$ by $H'\resto_{B'}=P$ and
$H'(x)=H(x) + (1-t)a +t a'+(dt)b$.

For any $Q\in\pi_0\shom_{B'}(B,A)$, the orbit map
$h^r(A)\to\pi_0\shom_{B'}(B,A)$ is equal to our earlier bijection described
in Lemma~\ref{Main}, defined using $Q$ as a base point for
$\shom_{B'}(B,A)$.  Thus our action is simply transitive. 
Note that it is given by
$$(a\ast Q)(x)=Q(x)+a\,.$$

Let us turn our attention to $\xi_P$.
The fact that $\xi_P$ is well-defined and a homomorphism is easily
checked.

Define a homomorphism
\begin{align}\label{need}
\pi_1\shom_C(B',A)&\longrightarrow h_1\big(Z^{r+1}
N(A\otimes\Omega\lcom)\big) \\*
h&\longmapsto h(dx)-P(dx)\,.\nonumber
\end{align}
To see that this is well-defined, let $h:B'\to
A\otimes\Omega_1$ and $h':B'\to A\otimes\Omega_1$ represent loops in
$\pi_1\shom_C(B',A)$. Let $\Psi:B'\to A\otimes\Omega_2$ be a homotopy
between them, i.e., let us assume that $\del_0\Psi=P$, $\del_1\Psi=h$
and $\del_2\Psi=h'$.  Then $\psi=\Psi(dx)-P(dx)$ is an element of
$Z^{r+1}(A\otimes\Omega_2)$ such that $\del_0\psi=0$,
$\del_1\psi=h(dx)-P(dx)$ and $\del_2\psi=h'(dx)-P(dx)$.
Then
$\tilde{\psi}(t_1,t_2)=\psi(t_1,t_2)-\big(h(dx)-P(dx)\big)(t_1+t_2)$
is an element $\tilde{\psi}\in Z^{r+1}N(A\otimes\Omega_2)$ such that
$\tilde\del(\tilde{\psi}) = \del_2(\tilde{\psi})=
\big(h'(dx)-P(dx)\big)-\big(h(dx)-P(dx)\big)$, showing that
$h(dx)-P(dx)$ and $h'(dx)-P(dx)$ are equal in $h_1\big(Z^{r+1}
N(A\otimes\Omega\lcom)\big)$.

We can show that (\ref{need}) is a homomorphism by a similar argument:
let $h$ and $h'$ be as before, except for not necessarily
homotopic. Let $\Psi:B'\to A\otimes\Omega_2$ be any 2-simplex such
that $\del_0\Psi=h$ and $\del_2\Psi=h'$.  Then $\del_1\Psi$ represents
$h\ast h'$, so without loss of generality we may set $h\ast
h'=\del_1\Psi$. Let us abbreviate notation to $\alpha=h(dx)-P(dx)$,
$\beta=h'(dx)-P(dx)$ and $\gamma=(h\ast h')(dx)-P(dx)$.  We need to
show that $\alpha+\beta-\gamma$ represents zero in $h_1\big(Z^{r+1}
N(A\otimes\Omega\lcom)\big)$. We again use the notation
$\psi=\Psi(dx)-P(dx)$. Then we have $\del_0\psi=\alpha$,
$\del_1\psi=\gamma$ and $\del_2\psi=\beta$.  This time we set
$\tilde{\psi}(t_1,t_2)=
\psi(t_1,t_2)+(\alpha-\gamma)(t_1+t_2)-\alpha(t_2)$. Then
$\tilde{\psi}\in Z^{r+1}N(A\otimes\Omega_2)$ and
$\tilde{\del}(\tilde{\psi})=\alpha+\beta-\gamma$, showing that
$\alpha+\beta-\gamma$ is a boundary, as required.

Now consider the diagram
$$\xymatrix{
h^{-1}\Deru_C(B',A)\dto_{\Xi_1}\rrto^{\xi_P} && h^r(A)
\dto^{\text{\upshape can.\ isom.\ (\ref{th.es})}}\\
\pi_1\shom_C(B',A)\rrto^{(\ref{need})} && h_1\big(Z^{r+1}
N(A\otimes\Omega\lcom)\big)}$$
It is easily checked that the diagram is commutative. We will
eliminate $h_1\big(Z^{r+1}N(A\otimes\Omega\lcom)\big)$ from this
diagram and replace it by
\begin{equation}\label{repl}
\vcenter{\xymatrix@R=10pt{
h^{-1}\Deru_C(B',A)\ddto_{\Xi_1}\drrto^{\xi_P} &&\\
&& h^r(A)\\
\pi_1\shom_C(B',A)\urrto_{\rho} &&}}
\end{equation}
where we have given the composition of (\ref{need}) with the inverse
of (\ref{th.es}) the name $\rho$. Let us emphasize again that the two
sloped maps in this diagram are homomorphisms, whereas $\Xi_1$ is
probably not a homomorphism, but it is bijective, because of
Theorem~\ref{Bij}. 

Note that the action of $\pi_1\shom_C(B',A)$ on $\pi_0\shom_{B'}(B,A)$
defined via $\rho$ is equal to the monodromy action indicated
in~(\ref{tail}).  To see this, let $h\in\pi_1\shom_C(B',A)$ and
$Q\in\pi_0\shom_{B'}(B,A)$. We need to prove that $h\ast Q=\rho(h)\ast
Q$. We may take $Q$ as base point for  $\pi_0\shom_{B'}(B,A)$. Then we
obtain the boundary map
$\del:\pi_1\shom_C(B',A)\to\pi_0\shom_{B'}(B,A)$ and we have $h\ast
Q=\del h$. Moreover, because $\Xi_1$ is bijective, we may assume that
$h=P+\omega_1 D=P+dtD$, for $D\in h^{-1}\Deru_C(B',A)$. We have seen
in the proof of Proposition~\ref{boundary} that then $\del h (x)
=Q(x)+\omega_0D(dx)=Q(x)+D(dx)$. Thus we have $(h\ast
Q)(x)=Q(x)+D(dx)$. On the other hand, by the commutativity
of~(\ref{repl}), we have $\rho(h)=D(dx)$. So by the definition of the
action of $h^r(A)$ on $\pi_0\shom_{B'}(B,A)$, we have
$\big(\rho(h)\ast Q\big)(x)=Q(x)+\rho(h)=Q(x)+D(dx)$. Since $h\ast Q$
and $\rho(h)\ast Q$ agree on $x$, they are, indeed, equal.

By the exactness properties of~(\ref{tail}), $h^r(A)$ acts
transitively on the fiber of $\pi_0\shom_C(B,A)\to\pi_0\shom_C(B',A)$
and the stabilizer is equal to the image of $\rho$.  Because of
Diagram~(\ref{repl}), and the bijectivity of $\Xi_1$, this image is
equal to the image of $\xi_P$.
\end{pf}

\begin{cor}\label{sim.tr}
The vector space $\cok \xi_P$ acts simply transitively on the fiber of 
$\pi_0\shom_C(B,A)\to\pi_0\shom_C(B',A)$
over $P$.  
\end{cor}

Both in the proposition and the corollary, there is no reason why the
fiber in question should be non-empty.  We will address the question
of when this fiber is non-empty next.

Let us fix $C\to B'\to B$ and $A$ as above, but let us forget about
$P$.

\begin{prop}\label{ext}
Define a map
\begin{align*}
\pi_0\shom_C(B',A)&\longrightarrow h^{r+1}(A)\\*
h&\longmapsto h(dx)\,.
\end{align*}
Then the sequence
$$\pi_0\shom_C(B,A)\longrightarrow \pi_0\shom_C(B',A)\longrightarrow
h^{r+1}(A)$$ 
is exact in the middle.
\end{prop}
\begin{pf}
Let $h:B'\to A$ represent an element of $\pi_0\shom_C(B',A)$. It lifts
to $\pi_0\shom_C(B,A)$ if and only if there exists a morphism of
differential graded algebras $h':B\to A$, such that $h'\resto_{B'}=h$,
because of the fact that $\shom_C(B,A)\to\shom_C(B',A)$ is a
fibration.

Suppose an extension $h'$ exists. Then we have
$h(dx)=h'(dx)=d\big(h'(x)\big)$, so that $h(dx)=0$ in $h^{r+1}(A)$. 

Conversely, assume that $h(dx)=da$, for some $a\in A^r$. Then we may
define an extension $h'$ of $h$ to $B$ to be the unique extension of
$h$ as a morphism of graded $k$-algebras satisfying $h'(x)=a$. This
respects the differential, because $d\big(h'(x)\big)=da=h(dx)=h'(dx)$.
\end{pf}

\begin{cor}
The map 
$$\pi_0\shom_C(B,A)\longrightarrow \ker\big(\pi_0\shom_C(B',A)\to
h^{r+1}(A)\big)$$
is surjective. For all $P$ in the kernel, the fiber over $P$ is a
principal homogeneous $cok\xi_P$-space.
\end{cor}

\begin{rmk}
We may also summarize the results of Propositions~\ref{xi}
and~\ref{ext} in 
the exact sequence
\begin{multline*}
\xymatrix{
h^{-1}\Deru_C(B',A)\rto&h^r(A)\ar@{o}[r]&}\\
\xymatrix{\ar@{o}[r] &
\pi_0\shom_C(B,A)\rto & \pi_0\shom_C(B',A)\rto & h^{r+1}(A)\,.}
\end{multline*}
This means that there is an action of $h^r(A)$ on $\pi_0\shom_C(B,A)$,
such that the orbits are equal to the fibers of the fibration
$\pi_0\shom_C(B,A)\to\pi_0\shom_C(B',A)$.  The stabilizers depend only
on the orbit.  For the orbit which is equal to the fiber over $P:B'\to
A$, the stabilizer is equal to the image of $h^{-1}\Deru_C(B',A)$
(whose definition depends on $P$) in $h^r(A)$.  Finally, a point of
$\pi_0\shom_C(B',A)$ has a non-empty fiber over it, if and only if it
maps to zero in $h^{r+1}(A)$. 
\end{rmk}

\subsection{Applications to \'etale morphisms}

Let $C\to B$ be a morphism of quasi-finite resolving algebras and
$B\to k$ an augmentation. The augmentation provides is with a
canonical base point for $\shom(B,A)$ and $\shom(C,A)$, for all
differential graded algebras $A$.

\begin{prop}
Let $r>0$ be an integer.
If $B$ and $C$ are perfect resolving algebras then the following are
equivalent: 

(i) $C\to B$ is \'etale at the augmentation $B\to k$,

(ii) the map of pointed spaces $\shom(B,A)\to\shom(C,A)$ induces an
isomorphism on homotopy groups $\pi_\ell$, for all $\ell\geq1$, and
all (finite) resolving algebras $A$,

(iii)  $\pi_r\shom(B,A)\to\pi_r\shom(C,A)$ an isomorphism, for all
(finite) resolving algebras $A$.

If, on the other hand, $C\to B$ is a perfect resolving morphism, then
$\shom(B,A)\to\shom(C,A)$ is a fibration, for every (finite) resolving
algebra $A$, and the following are equivalent:

(i) $C\to B$ is \'etale at the augmentation $B\to k$,

(ii) the fiber $\shom_C(B,A)$ is acyclic for all (finite) resolving
algebras $A$, 

(iii) $\pi_r\shom_C(B,A)=0$, for all  (finite) resolving algebras
$A$. 
\end{prop}
\begin{pf}
Simply combine Corollary~\ref{greatc} with Proposition~\ref{eac}.
\end{pf}

Let us now forget the augmentation $B\to k$. 

\begin{prop}\label{etahom}
Let $r>0$ be an integer.
If $B$ and $C$ are perfect resolving algebras then the following are
equivalent: 

(i) $C\to B$ is \'etale,

(ii) $\shom(B,A)\to\shom(C,A)$ induces isomorphisms on homotopy groups
$\pi_\ell$, for all $\ell\geq1$, and all (finite) resolving algebras
$A$ and all base points $B\to A$ 
for $\shom(B,A)$,

(iii)  $\pi_r\shom(B,A)\to\pi_r\shom(C,A)$ an isomorphism, for all
$B\to A$ as in~(ii).

If, on the other hand, $C\to B$ is a perfect resolving morphism, then
the following are equivalent: 

(i) $C\to B$ is \'etale,

(ii) the fiber $\shom_C(B,A)$ is acyclic for all (finite) resolving
algebras $A$  and all base points $B\to A$,

(iii) $\pi_r\shom_C(B,A)=0$, for all  $B\to A$ as in~(ii).
\end{prop}
\begin{pf}
This time,  combine Corollary~\ref{greatc} with Proposition~\ref{mea}.
\end{pf}

\Section{Finite resolutions}
\label{sec.fibered}

Our goal in this section is to prove that any morphism of finite
resolving algebras admits a finite resolution.  This is a significant
strengthening of Proposition~\ref{res.ex} in the finite case. As an
application, we can prove that the derived tensor product of finite
resolving algebras may be represented by a finite resolving algebra.

The existence of resolutions is reduced to the existence of `cylinder
objects' by a formal (and standard) argument. A cylinder object is
nothing other than a resolution of the diagonal $A\otimes A\to A$.

Let $A$ be a finite resolving algebra. Let $(x_i)_{i=1,\ldots,n}$ be a
basis for 
$A$ with the property that $i\leq j$ implies $\deg x_i\geq\deg
x_j$. For all $i=0,\ldots,n$ let $A_{(i)}$ denote the subalgebra of
$A$ generated by $x_1,\ldots,x_i$. Because of our assumption on the
degrees of the $x_i$, we have that $A_{(i)}$ is a differential graded
subalgebra of $A$, which is itself a finite resolving algebra with
basis $(x_1,\ldots,x_i)$. Notation: $dx_i=f_i(x)$. Note that
$f_i(x)=f_i(x_1,\ldots,x_{i-1}) \in A_{(i-1)}$.

The differential graded algebra $A\otimes A$ is a finite resolving
algebra with basis consisting of $y_i=x_i\otimes 1$ and $z_i=1\otimes
x_i$.

Let $\xi_i$, for $i=1,\ldots, n$, be a formal variable of degree $\deg
\xi_i=\deg x_i-1$, and let us consider the symmetric algebra $k[\xi]$
on the graded vector space with basis $(\xi_i)$. We also consider the
graded algebra $A\otimes
A\otimes k[\xi]=A\otimes A[\xi]$, with its subalgebra $A\otimes A$.

\begin{prop}\label{res.diag}
There is a way to extend the differential from the subalgebra
$A\otimes A$ to all of $A\otimes A[\xi]$, and to extend the diagonal
morphism $\Delta:A\otimes A\to A$ to all of $A\otimes A[\xi]$ in such
a way that $A\otimes A[\xi]$ becomes a differential graded algebra and
$A\otimes A[\xi]\to A$ a quasi-isomorphism of differential graded
algebras.
\end{prop}
In other words, $A\otimes A[\xi]$ will be a finite
resolution of the diagonal $A\otimes A\to A$.
\begin{pf}
The proof is by induction on $n$, the case $n=0$ serving as trivial
base case.  We assume that the proposition has been proved in the case
of $n$ generators, and we wish to prove that it is also true for the
case of $n+1$ generators.  Thus $d\xi_i$ and $\Delta(\xi_i)$, for
$i=1,\ldots,n$, have
already been found, making $\Delta:A_{(n)}\otimes
A_{(n)}[\xi_1,\ldots,\xi_n]\to A_{(n)}$ a quasi-isomorphism.  
We have $A=A_{(n+1)}$ and we need to
find suitable values for $d\xi_{n+1}$ and $\Delta(\xi_{n+1})$. Let
$r=\deg x_{n+1}$.  

\noindent{\sc Claim}\@. There exist $h(y,z,\xi)\in
A_{(n)}\otimes A_{(n)}[\xi_1,\ldots,\xi_n]$ of degree $r$ and $g(x)\in
A_{(n)}^{r-1}$ 
such that  

(i) $dh(y,z,\xi)=f_{n+1}(y)-f_{n+1}(z)$,

(ii) $\Delta(h)=h(x,x,\Delta\xi)=dg(x)$.

\noindent To prove that claim, we start by observing that $f_{n+1}(x)$
is a cocycle, because it is equal to $d x_{n+1}$ and $d^2=0$. Hence
$f_{n+1}(y)-f_{n+1}(z)$ is a cocycle, too. This cocycle maps to $0$
under $\Delta$, and so by the injectivity of the induced map
$$h^{r+1}\big(A_{(n)}\otimes A_{(n)}[\xi_1,\ldots,\xi_n]\big)\to
h^{r+1}(A_{(n)})\,,$$ 
we can find $g(y,z,\xi)\in A_{(n)}\otimes A_{(n)}[\xi_1,\ldots,\xi_n]$
of degree $r$, such that 
$$dg(y,z,\xi)=f_{n+1}(y)-f_{n+1}(z)\,.$$
Now we observe that $g(x,x,\Delta\xi)\in A$ is a cocycle, and so by the
surjectivity of 
$$h^{r}\big(A_{(n)}\otimes A_{(n)}[\xi_1,\ldots,\xi_n]\big)\to
h^{r}(A_{(n)})\,,$$ 
there exists a cocycle $\tilde{g}(y,z,\xi)\in A_{(n)}\otimes
A_{(n)}[\xi_1,\ldots,\xi_n]$, such that 
$$\tilde{g}(x,x,\Delta\xi)=g(x,x,\Delta\xi)\in h^r(A_{(n)})\,.$$
Then $h=g-\tilde{g}$ satisfies Condition~(i) of the claim. Moreover,
$\Delta(h)$ is a coboundary, so that we may choose $g(x)\in A_{(n)}$ such
that Condition~(ii) of the claim is also true. 

Now we set
$$d\xi_{n+1}=z_{n+1}-y_{n+1}+h(y,z,\xi)\,.$$
This turns $A\otimes A[\xi]$ into a differential
graded algebra, because $z_{n+1}-y_{n+1}+h$ is a cocycle:
$$d\big(z_{n+1}-y_{n+1}+h(y,z,\xi)\big)=
f_{n+1}(z)-f_{n+1}(y)+f_{n+1}(y)-f_{n+1}(z)=0\,.$$ 
Now we define
$$\Delta(\xi_{n+1})=g(x)\,.$$
This defines a morphism of differential graded algebras
$\Delta:A\otimes A[\xi]\to A$,
because 
$$\Delta(d\xi_{n+1})=\Delta(z_{n+1}-y_{n+1}+h)=\Delta(h)=
dg(x)=d(\Delta\xi_{n+1})\,.$$  
We claim that $\Delta$ is a quasi-isomorphism. We will prove this
using the criterion of Corollary~\ref{qiscondition}.

Let us first check that $h^0(\Delta)$ is an isomorphism.  This is
trivial if $r\leq-2$, because in that case $h^0$ is not affected when
passing from $n$ to $n+1$.  Let us consider the case $r=0$. Then
$h^0(A_{(n)})=A_{(n)}$ is a polynomial ring in $x_1,\ldots,x_n$ and
$h^0(A_{(n+1)})=A_{(n+1)}$ is a polynomial ring in
$x_1,\ldots,x_{n+1}$. On the other hand, $h^0(A_{(n)}\otimes
A_{(n)}[\xi_1,\ldots,\xi_n])$ is a quotient of the polynomial ring
$k[y_1,z_1,\ldots,y_n,z_n]$ by $n$ relations
$h_1(y,z),\ldots,h_n(y,z)$ and by the induction hypothesis we know
that $y_i\mapsto x_i$ and $z_i\mapsto 
x_i$ induces an isomorphism 
$$k[y_1,z_1,\ldots,y_n,z_n]/(h_1,\ldots,h_n)
\longrightarrow k[x_1,\ldots,x_n]\,.$$
This means that as ideals we have
$(h_1,\ldots,h_n)=(z_1-y_1,\ldots,z_n-y_n)$. Now when adjoining
$y_{n+1}$, $z_{n+1}$ and $\xi_{n+1}$ to $A_{(n)}\otimes
A_{(n)}[\xi_1,\ldots,\xi_n]$, we see that 
$$h^0(A\otimes
A[\xi])=k[y_1,z_1,\ldots,y_{n+1},z_{n+1}]/
(z_1-y_1,\ldots,z_n-y_n,z_{n+1}-y_{n+1}+h)\,,$$   
where $h(y_1,z_1,\ldots,y_{n},z_{n})$ satisfies
$h(x_1,x_1,\ldots,x_n,x_n)=0$, and hence
$h(y_1,z_1,\ldots,y_{n},z_{n})\in(z_1-y_1,\ldots,z_n-y_n)$. Therefore
we have
\begin{align*}
h^0(A\otimes
A[\xi])&=k[y_1,z_1,\ldots,y_{n+1},z_{n+1}]/
(z_1-y_1,\ldots,z_n-y_n,z_{n+1}-y_{n+1}+h)\\
&=k[y_1,z_1,\ldots,y_{n+1},z_{n+1}]/
(z_1-y_1,\ldots,z_{n+1}-y_{n+1})\\
&=k[x_1,\ldots,x_n]\\
&=h^0(A)\,.
\end{align*}
In the case $r=-1$ we have that
\begin{equation}\label{rmin}
h^0(A)= k[x_1,\ldots,x_m]/\big(f_{m+1}(x),\ldots,
f_{n+1}(x)\big)\,,
\end{equation}
where $m\leq n$ is the total number of degree $0$ generators among the
$x_i$. On the other hand, by the above results in the $r=0$ case, we
have 
\begin{multline*}
h^0(A\otimes A[\xi])=k[y_1,z_1,\ldots,y_{m},z_{m}]/\\
/\big(y_1-z_1,\ldots,y_m-z_m,f_{m+1}(y),f_{m+1}(z),\ldots,
f_{n+1}(y),f_{n+1}(z)\big)\,,
\end{multline*}
which is clearly equal to (\ref{rmin}).

Finally, we need to check that the cotangent complex of $\Delta$ is
acyclic.  Let us abbreviate $B=A\otimes A[\xi]$ and
$B_{(n)}=A_{(n)}\otimes A_{(n)}[\xi_1,\ldots,x_n]$.

Note that $\Omega_{B/B_{(n)}}\otimes_B h^0(B)$ is freely generated over
$h^0(B)$ by $\du y_{n+1}$, $\du z_{n+1}$ in degree $r$ and $\du
\xi_{n+1}$ in degree $r-1$. Moreover, 
$$d(\du \xi_{n+1})=\du(d\xi_{n+1})= \du z_{n+1}-\du y_{n+1} + \du
h(y,z,\xi)=\du z_{n+1}-\du y_{n+1}\,,$$
because $h(y,z,\xi)\in B_{(n)}$. 

On the other hand, $\Omega_{A/A_{(n)}}\otimes_A h^0(A)$ is freely
generated over $h^0(A)=h^0(B)$ in degree $r$ by $\du x_{n+1}$. 

The map $\du\Delta$ maps $\du y_{n+1}$ and $\du z_{n+1}$ to
$\du x_{n+1}$ and it maps $\du \xi_{n+1}$ to $\du g(x)=0$, because
$g(x)=A_{(n)}$. 

Thus it is clear that $\du\Delta:\Omega_{B/B_{(n)}}\otimes_B
h^0(A)\to \Omega_{A/A_{(n)}}\otimes_A h^0(A)$ is a quasi-isomorphism,
which is all that we needed to prove.
\end{pf}

\begin{cor}\label{twomh}
The two morphisms $A\to A\otimes A[\xi]$ given by $a\mapsto a\otimes1$
and $a\mapsto 1\otimes a$ are homotopic. (They are also
quasi-isomorphisms.)
\end{cor}
\begin{pf}
Both of these morphisms are 2-inverses of the quasi-isomorphism
$A\otimes A[\xi]\to A$, constructed in the proposition. Use
Corollary~\ref{qiseq}.
\end{pf}

\begin{them}\label{fin.res.ex}
Let $A\to B$ be a morphism of finite resolving algebras. Then there
exists a resolution $A\to B'\to B$, where $A\to B'$ is a finite
resolving morphism.
\end{them}
\begin{pf}
See the proof of Corollary~8.3 in Section~II of~\cite{simhomthe}.  In
our case, the resolution of $A\to B$ is given by 
\begin{align*}
A&\longrightarrow B\otimes A[\xi] 
\stackrel{B\otimes\sigma}{\longrightarrow} B\\*
a & \longmapsto 1\otimes a  
\end{align*}
where $\sigma:A\otimes A[\xi]\to A$ is any resolution of the diagonal.
\end{pf}

\begin{cor}
Let $f:A\to B$ and $g:A\to C$ be morphisms of finite resolving
algebras. Then there exists a finite resolving algebra $R$ and a
homotopy commutative square 
\begin{equation}\label{wfpp}
\vcenter{\xymatrix{
R\drtwocell\omit{^} & C\lto\\
B\uto & A\lto^{f}\uto_{g}}}
\end{equation}
such that, whenever we resolve $f$ as $A\to B'\to B$, then there
exists a factorization of (\ref{wfpp}) as
$$\xymatrix{
R\drtwocell\omit{^}& B'\otimes_AC\lto & C \lllowertwocell{\omit}\lto\\ 
B\uto & B'\uto\lto & A\lto\uto }$$
where $B'\otimes_A C\to R$ is a quasi-isomorphism.
\end{cor}
\begin{pf}
This is a formal consequence of Theorem~\ref{fin.res.ex}.  More
explicitly, we can construct $R$ as follows:
let $A$ be as above, generated by $x_1,\ldots,x_n$, with
$dx_i=f_i(x)$. Choose for every $i$ an element $h_i(y,z,\xi)\in A\otimes
A[\xi]$ as in the inductive proof of Proposition~\ref{res.diag}. In
particular, $dh_i(y,z,\xi)=f_i(y)-f_i(z)$.
Denote the images of the $x_i$ in $B$ by $b_i$
and in $C$ by $c_i$. Consider  the differential graded algebra 
$$R=B\otimes C[\xi]$$
with differential defined by $d\xi_i=c_i-b_i+h_i(b,c,\xi)$.
Any homotopy between the two morphisms $A\to A\otimes A[\xi]$ (see
Corollary~\ref{twomh}), composed with $A\otimes A[\xi]\to B\otimes
C[\xi]$ gives rise to a homotopy commutative diagram
\begin{equation*}\label{wfp}
\vcenter{\xymatrix{
B\otimes C[\xi] & C\lto\\
B\uto & A\lto^f\uto_{ g}\ultwocell\omit}}
\end{equation*}
as required.
\end{pf}

\begin{scholum}
If $a$, $b$ and $c$ denote the numbers of elements of finite bases for
$A$, $B$ and $C$, respectively, then there exists a basis for $R$ with
$a+b+c$ elements.
\end{scholum}

\begin{numrmk}
The results of this section also apply to perfect resolving algebras
and perfect resolutions, in place of finite resolving algebras and
finite resolutions.  The proofs are easier and require only the use of
Lemma~\ref{cot.cx}~(iii).
\end{numrmk}

%\nocite{*}

%\bibliographystyle{plain}
%\bibliography{basic}    

\end{document}